\documentclass[12pt, a4paper, reqno]{amsart}
\usepackage{amssymb}
\usepackage{amsthm}
\usepackage[dvipsnames]{xcolor}
\usepackage{tikz}
\usetikzlibrary{arrows,positioning} 
\usepackage{mathrsfs}
\usepackage{bbm}
\usepackage{ stmaryrd }
\usepackage[normalem]{ulem}
\usepackage{framed}

\usepackage{hyperref}
\hypersetup{
	colorlinks=true,
	linkcolor=blue,
	filecolor=magenta,  
	urlcolor=cyan,
}

\usepackage{hyperref}
\usepackage[capitalise]{cleveref}

\usepackage{amsfonts}
\usepackage{amsthm}
\usepackage{amsmath}
\usepackage{amscd}
\usepackage[latin2]{inputenc}
\usepackage{t1enc}
\usepackage[mathscr]{eucal}
\usepackage{indentfirst}
\usepackage{graphicx}
\usepackage{graphics}
\usepackage{epic}
\numberwithin{equation}{section}
\usepackage[margin=2.9cm]{geometry}
\usepackage{epstopdf} 
\usepackage{subfig}

\usepackage{multirow}
\usepackage{makecell}

\def\ker{\operatorname{ker}}

\def\im{\operatorname{Im}}

\def\dim{\operatorname{dim}}

\def\ad{\operatorname{ad}}

\def\Hom{\operatorname{Hom}}

\def\C{\mathbb{C}}

\def\N{\mathbb{N}}
\def\Z{\mathbb{Z}}
\def\G{\mathbb{G}}

\def\AA{\mathcal{A}}

\def\BB{\mathcal{B}}

\def\NN{\mathcal{N}}

\def\SS{\mathcal{S}}
\def\CC{\mathcal{C}}

\def\II{\mathcal{I}}

\def\b{\mathfrak{b}}

\def\p{\mathfrak{p}}
\def\q{\mathfrak{q}}
\def\g{\mathfrak{g}}

\def\t{\mathfrak{t}}
\def\h{\mathfrak{h}}

\def\k{\mathfrak{k}}
\def\l{\mathfrak{l}}
\def\s{\mathfrak{s}}
\def\o{\mathfrak{o}}

\def\ol{\overline}

\def\Gr{\text{Gr}}

\def\Rep{\operatorname{Rep}}

\def\sub{\subseteq}

\def\xto{\xrightarrow}

\newcommand{\SVec}{\mathtt{sVec}}

\newtheorem{thm}{Theorem}[section]
\newtheorem{cor}[thm]{Corollary}
\newtheorem{lemma}[thm]{Lemma}
\newtheorem{prop}[thm]{Proposition}

\theoremstyle{definition}
\newtheorem{definition}[thm]{Definition}
\theoremstyle{remark}
\newtheorem{remark}[thm]{Remark}
\newtheorem{example}[thm]{Example}

\numberwithin{equation}{section}

\def\quotient#1#2{%
	\raise1ex\hbox{$#1$}\Big/\lower1ex\hbox{$#2$}%
}

\newcommand{\InnaA}[1]{#1}
\newcommand{\InnaB}[1]{#1}

\newcommand{\comment}[1]{}

\usepackage{relsize}
\usepackage{fancyhdr}
\usepackage[all,cmtip]{xy}

\begin{document}
	
	\title{It takes two spectral sequences}
	
	\author{Inna Entova-Aizenbud, Vera Serganova, Alexander Sherman}

	\pagestyle{plain}
	\maketitle

	{\centering\footnotesize Dedicated to the memory of Georgia Benkart.\par}

	\begin{abstract} We study the representation theory of the Lie superalgebra $\g\l(1|1)$, constructing two spectral sequences which eventually annihilate precisely the superdimension 0 indecomposable modules in the finite-dimensional category. The pages of these spectral sequences, along with their limits, define symmetric monoidal functors on $\Rep \g\l(1|1)$. These two spectral sequences are related by contragredient duality, and from their limits we construct explicit semisimplification functors, which we explicitly prove are isomorphic up to a twist. We use these tools to prove branching results for the restriction of simple modules over Kac-Moody and queer Lie superalgebras to $\g\l(1|1)$-subalgebras.  
	\end{abstract}
	
	\section{Introduction}

	\subsection{}

	Consider the general linear Lie superalgebra $\g\l(1|1)$ presented as matrices of the form 
	\[
	\begin{bmatrix}c+h/2 & x\\ y & c-h/2\end{bmatrix}.
	\]
	Given a $\g\l(1|1)$-module $V$ on which $h$ and $c$ act semisimply, if we take invariants nder $c$ we obtain a super vector space $V^c$ with a supercommuting action of $x$ and $y$, along with a semisimple action of $h$.  The supercommuting actions of $x$ and $y$ allow one to write down a natural double complex whose terms are all $M:=V^c$, and then consider the two spectral sequences associated to it.  Let us consider the one in which we take cohomology with respect to $x$ first.  All entries of this spectral sequence are the same at any given page, so picking out what happens at a fixed position, we obtain a sequence of semisimple $\C\langle h\rangle$-modules which we write as $DS_{x,y}^nM$.  Here $DS$ refers to Duflo-Serganova, in reference to the connection of our functors to the Duflo-Serganova functors in the representation theory of Lie superalgebras.
	
	As one might expect, the assignment $M\mapsto DS_{x,y}^nM$ is functorial; in fact it defines a symmetric monoidal functor to $\Rep\C\langle h,d_{1-2n}\rangle$, which is the category of \InnaA{finite-dimensional} representations of the Lie superalgebra with even generator $h$ acting semisimply, odd generator $d_{1-2n}$, and such that
	\[
	[h,d_{1-2n}]=(1-2n)d_{1-2n},  \ \ \ [d_{1-2n},d_{1-2n}]=0.
	\]
	Thus we have a sequence of symmetric monoidal functors, where $DS_{x,y}^1M=M_x:=DS_xM$, and $DS_{x,y}^2M=(M_x)_{y}:=DS_y(DS_xM)$.  Here $DS_x$ and $DS_y$ are Duflo-Serganova functors, see \cite{GHSS22}.  
	
	In general, we have $DS_{x,y}^{n+1}M$ is the cohomology of $d_{1-2n}$ on $DS_{x,y}^nM$.  These functors detect, to some degree, the indecomposable components of $M$ as a $\g\l(1|1)$-module.  In particular if $M$ is finite-dimensional, the spectral sequence converges, and we write the object it converges to as $DS^{\infty}_{x,y}M$, where $DS^{\infty}_{x,y}$ itself defines a  symmetric monoidal functor.  We of course may do the same procedure with the other spectral sequence, whose terms we write as $DS^n_{y,x}$, which will define symmetric monoidal functors to $\Rep\C\langle h,d_{2n-1}\rangle$.
	
	\subsection{Relations between functors} The functors $DS^\infty_{x,y}$ and $DS_{y,x}^\infty$ are useful on their own, but they also have precise, useful connections to the functor $DS_{x+y}$ and semisimplification functors on $\Rep_{\g_{\ol{0}}}\g\l(1|1)$, i.e. the category of $\g\l(1|1)$-modules with semisimple action of $h$ and $c$. We summarize these connections in the diagram below.
	
	Note that $\Rep_{\g_{\ol{0}}}\g\l(1|1)$ is an abelian rigid symmetric monoidal category, and (as we will show) its semisimplification is $\Rep\g_m\times\Rep\G_m$.  We construct two semisimplification functors as depicted, $DS_{x,y}^{ss}$ and $DS_{y,x}^{ss}$; we show that these two functors are equivalent up to an autoequivalence $\Phi_{\phi}$, which comes from twisting by the automorphism $\phi$ of $\C\langle h\rangle\times\operatorname{Lie}\G_m$ given by $(h,a)\mapsto(h+2a,a)$. 
	
	A little more on notation: $\Rep\g_m^{fil}$ denotes the category of $\Z$-filtered modules in $\Rep\g_m$, and $\operatorname{Fil}^{\Z\times\C}$ is a certain category of filtered super vector spaces (see \cref{ssec:DSx+y}). The functors $\Gr$ and $\Gr_{a_i}$ are functors given by taking an associated graded space.  
	
	Note that all smaller, enclosed triangles in the diagram are commutative in the sense that we have a natural isomorphism of symmetric monoidal functors.   Further, contragredient duality acts on the picture by reflecting about the central vertical axis
	
	
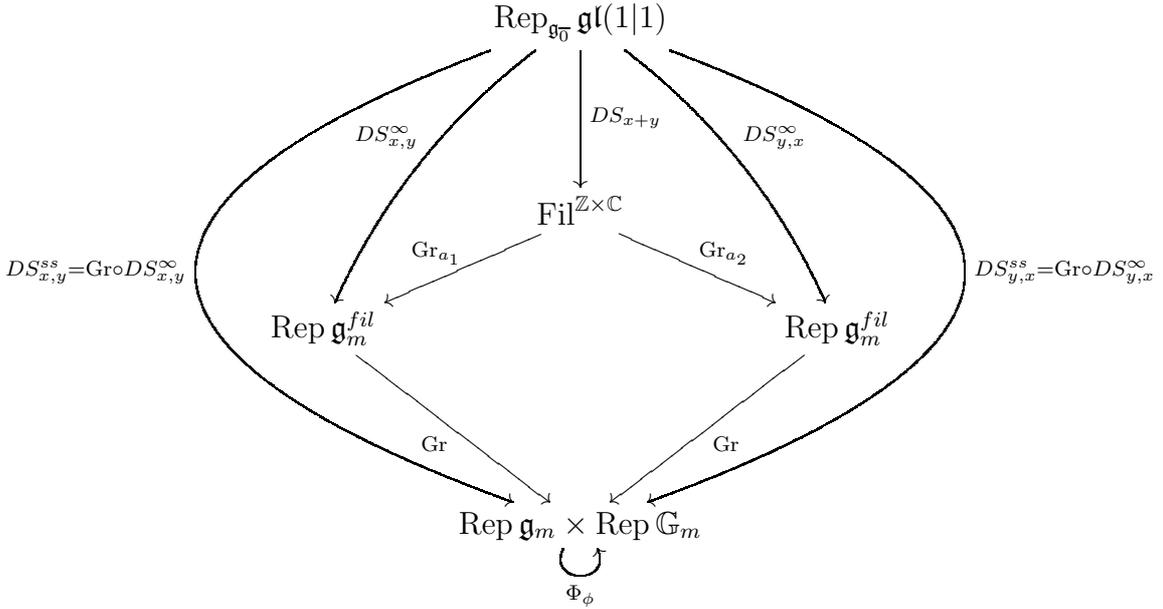
\[
	\xymatrix{ & \Rep_{\g_{\ol{0}}}\g\l(1|1) \ar[dd]^{DS_{x+y}}  \ar@/^12pc/[ddddd]^{DS_{y,x}^{ss}=\Gr\circ DS_{y,x}^\infty} \ar@/_12pc/[ddddd]_{DS_{x,y}^{ss}=\Gr\circ DS_{x,y}^{\infty}} \ar@/_1pc/[dddl]_{DS_{x,y}^{\infty}} \ar@/^1pc/[dddr]^{DS_{y,x}^\infty} & \\ 
		&&\\
		& \operatorname{Fil}^{\Z\times \C}\ar[dr]^{\Gr_{a_2}} \ar[dl]_{\Gr_{a_1}} & \\
		\Rep\g_m^{fil}\ar[ddr]_{\operatorname{Gr}} & & \Rep\g_m^{fil} \ar[ddl]^{\operatorname{Gr}}\\
		&&\\
		& \Rep\g_m\times\Rep\G_m  \ar@(dl,dr)_{\Phi_{\phi}} &
	}
	\]
	
	\begin{remark}
		We note that the top half of the above diagram can be thought of morally as a consequence of the theory of spectral sequences.  This theory tells us that the cohomology of the total complex, which should correspond to $DS_{x+y}$ in the above picture, has two filtrations whose associated graded ($\Gr_{a_1}$ and $\Gr_{a_2}$ in the above) give the last pages of the  two spectral sequences ($DS_{x,y}^\infty$ and $DS_{y,x}^\infty$).
	\end{remark}
	
	Written in more explicit terms we have the following theorem summarizing the situation:
	\begin{thm}
		\mbox{}
		\begin{enumerate}
			\item The functors $DS^{\infty}_{x,y}$ and $DS_{y,x}^\infty$ define full, essentially surjective functors \linebreak $\Rep_{\g_{\ol{0}}}\g\l(1|1)\to\Rep\g_m^{fil}$ such that $DS_{x,y}^\infty\circ(-)^\vee\cong (-)^\vee\circ DS_{y,x}^\infty$, where $(-)^\vee$ denotes the contragredient duality functor.
			\item If we denote by $\Gr:\Rep\g_m^{fil}\to\Rep\g_m\times\Rep\G_m$ the functor of taking associated graded, then $DS_{x,y}^{ss}:=\Gr\circ DS^{\infty}_{x,y}$ and $DS^{ss}_{y,x}:=\Gr\circ DS^{\infty}_{y,x}$ define semisimplification functors, and we have an isomorphism of functors
			\[
			DS_{x,y}^{ss}\cong \Phi_{\phi}\circ DS_{y,x}^{ss}.
			\]
			\item The functor $DS_{x+y}$ defines a full, essentially surjective functor \linebreak $\Rep_{\g_{\ol{0}}}\g\l(1|1)\to\operatorname{Fil}^{\Z\times \C}$.  
			\item If we denote by $\Gr_{a_1},\Gr_{a_2}:\operatorname{Fil}^{\Z\times \C}\to\Rep\g_m^{fil}$ the two associated graded functors (defined in \cref{ssec:DSx+y}), we have natural isomorphisms of functors
			\[
			DS_{x,y}^{\infty}\cong \Gr_{a_1}\circ DS_{x+y},  \ \ \ \ DS_{y,x}^\infty\cong\Gr_{a_2}\circ DS_{x+y}.
			\]
		\end{enumerate}
	\end{thm}
	
	\subsection{Applications} Using the functors defined above and their connections to one another, we prove results about branching to a subalgebra $\g\l(1|1)$. 
	
	\begin{thm}
		Let $\g$ be a finite-type Kac-Moody Lie superalgebra or $\q(n)$, and consider a diagonal $\g\l(1|1)$-subalgebra of $\g$ (see Definition \ref{def_gl(11)_root_sub}).  Assume that $L$ is a simple finite-dimensional module, and if $\g=\q(n)$ also assume that its highest weight is integral: then the restriction to $\g\l(1|1)$ contains no non-projective indecomposable summands of superdimension zero. 
	\end{thm}
	
	We note that for $\q(n)$-modules of half-integral weight, non-projective indecomposable summands of superdimension 0 do appear; however for certain subalgebras $\g\l(1|1)$ we are able to compute which ones appear, and with what multiplcity (see \cref{ssec:half integral}).
	
	For a Lie superalgebra $\g$, write $\Rep_{\g_{\ol{0}}}^+(\g)$ for the Karoubian monoidal subcategory of $\Rep_{\g_{\ol{0}}}(\g)$ generated by all simple modules, where $\Rep_{\g_{\ol{0}}}(\g)$ is the symmetric monoidal category of $\g$-modules that are semisimple over $\g_{\ol{0}}$.
	
	\begin{thm}\label{thm intro branching}
		If $\g=\g\l(m|n)$ or $\o\s\p(2m|2n)$, and $\g\l(1|1)$ is a root subalgebra of $\g$, then for a simple finite-dimensional $\g$-module $L$ we have $\operatorname{Res}_{\g\l(1|1)}L\in\Rep_{\g_{\ol{0}}}^+(\g\l(1|1))$.
	\end{thm}
	See \cref{theorem branching} for the proof; we expect \cref{thm intro branching} to extend to $\o\s\p(2m+1|2n)$, however we do not yet see a way to prove this using our techniques.
	
	\subsection{Summary of article}  In Section 2 we recall facts about $\g\l(1|1)$-modules, including listing all indecomposables, which will be used throughout the article.  Section 3 defines the spectral sequences and gives the main results, namely that we obtain two sequences of symmetric monoidal functors $DS_{x,y}^n$ and $DS_{y,x}^n$ for all $n\in\N\cup\{\infty\}$.  In addition, contragredient duality is shown to interchange $DS_{x,y}^n$ and $DS_{y,x}^n$. Section 4 uses $DS_{x,y}^\infty$ and $DS_{y,x}^\infty$ to explicitly realize the semisimplification of $\Rep_{\g_{\ol{0}}}(\g\l(1|1))$, and proves directly the important relationship between the two semisimplification functors $DS_{x,y}^{ss}$ and $DS_{y,x}^{ss}$.  Further, the functor $DS_{x+y}$ is studied and its relationship to $DS_{x,y}^\infty$ and $DS_{y,x}^\infty$ is explained.  In Section 5 we give applications of our constructions, in particular proving results about branching simple $\g$-modules to $\g\l(1|1)$-submodules when $\g$ is finite-type Kac Moody or $\q(n)$.  Also for $\q(n)$ we explicitly compute the functors $DS_{x,y}^n$ on simple modules when $x,y$ are opposite root vectors.  Finally, the appendix gives all details of the spectral sequence.
	
	\subsection{Acknowledgements} The authors would like to thank Maria Gorelik, Thorsten Heidersdorf, and Vladimir Hinich for many helpful discussions. The authors were supported by the NSF-BSF grant 2019694.
	
	\section{Preliminaries and notation: \texorpdfstring{$\g\l(1|1)$}{gl(1|1)}, and \texorpdfstring{$\p\g\l(1|1)$}{pgl(1|1)}}\label{sec:prelim}

	\subsection{General notation}
	We work over field $\C$ of complex numbers; all categories and functors will be $\C$-linear. We will denote by $\SVec$ the category of $\C$-linear finite-dimensional vector superspaces and grading preserving maps. For a vector superspace $V$ we write $V=V_{\ol{0}}\oplus V_{\ol{1}}$ for its $\Z_2$-grading.  
	
	We will be considering quasireductive Lie superalgebras $\g$, i.e. those with $\g_{\ol{0}}$ reductive and $\g_{\ol{1}}$ a semisimple $\g_{\ol{0}}$ module.  We will write $\Rep_{\g_{\ol{0}}}(\g)$ for the category of finite-dimensional $\g$-modules which are semisimple over $\g_{\ol{0}}$. 
	
	We also denote by $\Rep^+_{\g_{\ol{0}}}(\g)$ the Karoubian symmetric monoidal subcategory of $\Rep_{\g_{\ol{0}}}(\g)$ generated by irreducible modules.
	
	\subsection{The category \texorpdfstring{$\Rep\g_m$}{T} }
	We denote by $\Rep\g_m$ the category of semisimple finite-dimensional $\C\langle h\rangle$-modules, and by $\Rep\g_m^{fil}$ the category of $\Z$-filtered, semisimple finite-dimensional $\C\langle h\rangle$-modules.
	
	Note that we do not require that $h$ act by integer eigenvalues.
	
	\subsection{Lie superalgebra \texorpdfstring{$\g\l(1|1)$}{gl(1|1)}}\label{ssec:gl_1_1} 
	
	Consider the Lie superalgebra $\g=\g\l(1|1)$ presented with basis $h,c,x,y$ such that 
	\[
	c\text{ is central}, \ \ \ [h,x]=x, \ \ \ [h,y]=-y, \ \ \ [x,y]=c.
	\]
	We see that $[\g,\g]=\langle c,x,y\rangle$, so that $\g/[\g,\g]\cong\C\langle h\rangle$.  Write $\mathbf{str}:\g\to \C$ for the character sending $h$ to $1$, so that for $r\in\C$ we have $r\mathbf{str}:\g\to\C$ with $r\mathbf{str}(h)=r$.  Then given a $\g$-module $V$, we write $V_r$ for its twist by the character $r\mathbf{str}$.
	
	The blocks of $\Rep_{\g_{\ol{0}}}(\g)$ either have a nontrivial action of $c$, in which case they consist of one simple projective module of dimension $(1|1)$ along with its parity shift, or $c$ acts trivially.  When $c$ acts by 0, we obtain the representations of $\p\g\l(1|1)=\g\l(1|1)/\langle c\rangle$.  Observe that we have an exact functor 
	\[
	\Rep_{\g_{\ol{0}}}\g\l(1|1)\to\Rep_{\g_{\ol{0}}}\p\g\l(1|1),  \ \ \ \ \ V\mapsto V^c,
	\]
	where $V^c$ denotes taking invariance under $c$.
	
	The blocks on which $c$ acts trivially are indexed by $\C/2\Z$, corresponding the possible eigenvalues of $h$ on even vectors, modulo $2\Z$.  The principal block corresponds to $0\in\C/2\Z$.  Factoring out by $2\Z$ corresponds to the fact that twisting by $2\mathbf{str}$ preserves blocks with trivial central character.

	\subsection{Lie superalgebra \texorpdfstring{$\p\g\l(1|1)$}{pgl(1|1)} and its indecomposable representations}\label{ssec:indecomps} 
	
	From our presentation of $\g\l(1|1)$ above, we have that $\p\g\l(1|1)$ has generators $h$ even, $x,y$ odd, satisfying:	
	\[
	[h,x]=x,\ \ \ \ [h,y]=-y,\ \ \ \ [x,y]=0.
	\]
	We present below the finite-dimensional indecomposable modules (up to parity shift and weight shift) of $\p\g\l(1|1)$ with semisimple action of $h$, which follows from \cite{R75} or \cite{J98}; note that as explained above these are also the finite-dimensional modules over $\g\l(1|1)$ with trivial action of $c$.  When we refer to the weights of the module, we mean the eigenvalues of $h$. 
	\begin{itemize}
		\item The projective indecomposable $P$ on the trivial module.
		\[
		\xymatrix{& \bullet\ar[dr]^x\ar[dl]_y  & \\ \bullet\ar[dr]_x &&\bullet\ar[dl]^y \\ & \bullet & }
		\]
		\item $X(n)$, $n\in\Z_{>0}$, a $'Z'$-module: dimension $(n|n)$, lowest weight is $-n+1/2$ and even, highest weight is $n-1/2$ and odd.
		\[
		\xymatrix{\bullet\ar[dr]^x &  & \bullet\ar[dr]^x \ar[dl]_y & & \hdots & \bullet\ar[dr]^x& \\
			& \bullet & & \bullet & \hdots & & \bullet  }
		\]
		\item $Y(n)$, $n\in\Z_{>0}$, a $'Z'$-module: dimension $(n|n)$, lowest weight is $-n+1/2$ and even, highest weight is $n-1/2$ and odd.
		\[
		\xymatrix{ & \bullet \ar[dl]_y & \hdots & & \bullet\ar[dr]^x \ar[dl]_y & & \bullet\ar[dl]_y\\
			\bullet & &  \hdots & \bullet & & \bullet }
		\]
		\item $W(n)$ for $n\geq0$, a $'W'$-module: dimension $(n+1|n)$, lowest weight is $-n$, highest weight is $n$, and both are even.
		\[
		\xymatrix{\bullet\ar[dr]^x &  & \bullet\ar[dr]^x \ar[dl]_y & & \hdots & \bullet\ar[dr]^x& & \bullet \ar[dl]_y\\
			& \bullet & & \bullet & \hdots & & \bullet & }
		\]
		\item $W(-n):=W(n)^*$ for $n\geq 0$, a $'W'$-module: dimension $(n+1|n)$, lowest weight is $-n$, highest weight is $n$, and both are even.
		\[
		\xymatrix{\bullet &  & \bullet & & \hdots & \bullet& & \bullet \\
			& \bullet\ar[ur]^x\ar[ul]_y & & \bullet\ar[ul]_y & \hdots & & \bullet\ar[ur]^x\ar[ul]_y & }
		\]
	\end{itemize} 
	
	Observe that $W(0)$ is isomorphic to the trivial module.  Further we can twist any of these above modules by $r\mathbf{str}$, or take their parity shift.  This will give all finite-dimensional indecomposables in $\Rep_{\g_{\ol{0}}}\p\g\l(1|1)$. 
	
	\subsection{Contragredient duality}\label{ssec:contra duality} 
	
	On $\g\l(1|1)$ we have a Chevalley automorphism \linebreak $\sigma=\sigma_{\g\l(1|1)}$ given by minus supertranspose, i.e.
	\[
	\sigma(h)=-h, \ \ \sigma(c)=-c, \ \ \sigma(x)=-y, \ \ \sigma(y)=x.
	\]
	Further, $\sigma$ descends to an automorphism of $\p\g\l(1|1)$.  This automorphism induces contragredient duality $V\mapsto V^\vee:=(V^*)^{\sigma_{\g\l(1|1)}}$.  It acts on indecomposables in the following way (up to isomorphism):
	\[
	P\mapsto P, \ \ \ \ X(n)\mapsto Y(n), \ \ \ \ Y(n)\mapsto X(n), \ \ \ \ W(n)\mapsto W(-n).
	\]	
	For any finite-dimensional $\g\l(1|1)$-module we have a natural isomorphism
	\[
	V\mapsto (V^\vee)^\vee,
	\]
	given by
	\[
	v\mapsto(-1)^{\ol{v}}v.
	\]
	\subsection{Lemma on properties of certain morphisms}
	
	\begin{lemma}\label{lemma morphism} Let $r,s\in\C$.
		\begin{enumerate}
			\item Suppose that $m<n\leq0$.  Then any map $W(m)_r\to W(n)_s$ must contain the socle of $W(m)_r$ in its kernel.	
			\item Suppose that $0\leq m<n$.  Then the image of any map $W(m)_r\to W(n)_s$ is contained in the radical (equivalently, the socle) of $W(n)_s$.
			\item For $m\leq 0$, $n>0$, any maps $X(n)_r\to W(m)_s$ or $Y(n)_r\to W(m)_s$ must vanish on the socle.
			\item For $m\geq0$, $n>0$, any maps $W(m)_s\to X(n)_r$ or $W(m)_s\to Y(n)_r$ must have image contained in the radical (equivalently, the socle).
		\end{enumerate}
	\end{lemma}
	\begin{proof}
		Notice that $(1)\iff(2)$ via contragredient duality.  To prove (1), assume without loss of generality that $r=0$, and $s\in\Z$ and let $\phi:W(m)\to W(n)_s$ be a morphism with $m<n\leq 0$.  Write $v_{m},\dots,v_{-m}$ for a weight basis of $W(m)$.  Then necessarily one of $v_{m}$ or $v_{-m}$ must lie in the kernel of $\phi$; let us assume that $\phi(v_{m})=0$, with the other case being a similar argument.  Then necessarily $y\phi(v_{m+1})=0$; we see that from the structure of $W(n)_s$ this forces $\phi(v_{m+1})$ to lie in the socle of $W(n)_s$, which implies that $\phi(v_{m+2})=x\phi(v_{m+1})=0$.  Thus we may continue inductively to obtain that $\phi(v_{m+2i})=0$ for all $i$, proving (1).
		
		Again by contragredient duality we have that $(3)\iff(4)$, so we show (3).  Further, since $X(n)^{\sigma}\cong Y(n)$ and $W(n)^{\sigma}\cong W(n)$, it suffices to deal with $X(n)$.  Assume WLOG that $r=0$ and $s\in\Z$, and let $\phi:X(n)\to W(m)_s$ be a morphism.  Write $v_{-n+1/2},\dots,v_{n-1/2}$ for a weight basis of $X(n)$; then since $y\phi(v_{-n+1/2})=0$, we must have $\phi(v_{-n+1/2})$ lies in the socle, implying that $\phi(v_{-n+3/2})=x\phi(v_{-n+1/2})=0$.  This in turn implies that $y\phi(v_{-n+5/2})=0$, and we may again continue the argument inductively.
	\end{proof}
	
	\subsection{Duflo-Serganova functors}  
	\subsubsection{Definition}
	
	The Duflo-Serganova functors were introduced in \cite{DS05}; see the recent survey \cite{GHSS22} for more on the role of this functor in the representation theory of Lie superalgebras. 
	
	We consider a slight extension of this functor. For a Lie superalgebra $\g$, set $x^2:=\frac{1}{2}[x,x]$ for $x\in\g_{\ol{1}}$, and define
	\[
	\g^{hom}_{\ol{1}}:=\{x\in\g:\ad(x^2)\text{ is semisimple}\}.
	\]
	Let $x\in \g^{hom}_{\ol{1}}$, and take $M$ to be a $\g$-module on which $x^2$ acts semisimply. Then we define
	\[
	DS_xM:=M_x=\frac{\ker(x:M^{x^2})}{\im(x:M^{x^2})}.  
	\]
	In particular $\g_x$ will be a Lie superalgebra, and $M_x$ will have the natural structure of a $\g_x$-module.  
	
	\InnaA{This defines a functor $$DS_x: \Rep(\g) \to \Rep(\g_x)$$ which is symmetric monoidal and $\C$-linear, but not exact: given a short exact sequence of $\g$-modules $$0 \to M' \to M \to M'' \to 0$$ we obtain an exact sequence $$	M'_x \to M_x \to M''_x$$ which is not necessarily exact on either side (see \cite{GHSS22}). We call this property of the Duflo-Serganova functors "middle exactness". 
	}
	\subsubsection{\texorpdfstring{$DS$}{DS} functors on \texorpdfstring{$\g\l(1|1)$}{gl(1|1)}-modules}  For $\g=\g\l(1|1)$ we have $\g^{hom}_{\ol{1}}=\g_{\ol{1}}$.  The following lemma is straightforward.
	\begin{lemma}\label{lemma DS on indecomps}
		We have the following:
		\begin{enumerate}
			\item $DS_uQ=0$ for all non-zero $u\in\g\l(1|1)_{\ol{1}}$ and all projective modules $Q$	.
			\item $DS_xX(n)=0$, $DS_yX(n)$ is $(1|1)$-dimensional.
			\item $DS_yY(n)=0$, $DS_xY(n)$ is $(1|1)$-dimensional.
			\item $DS_uW(n)$ is $(1|0)$-dimensional for all $n\in\Z$ and all non-zero $u\in\p\g\l(1|1)_{\ol{1}}$.
		\end{enumerate}
	\end{lemma}

	\section{A spectral sequence}\label{sec:spectral_seq}
	
	\subsection{The setting}\label{ssec:spectral_seq_intro}
	Let $M$ be an $\p\g\l(1|1)$-module, of any dimension. Then we have a double complex
	\[
	\xymatrix{
		& \vdots \ar[d] & \vdots \ar[d] & \vdots \ar[d] & \\
		\hdots \ar[r] & M \ar[r]^x\ar[d]^y &  M\ar[d]^y\ar[r]^x & M\ar[r]\ar[d]^y & \hdots \\
		\hdots \ar[r] &  M \ar[r]^x\ar[d]^y & M\ar[d]^y\ar[r]^x &  M\ar[r]\ar[d]^y & \hdots \\
		\hdots \ar[r] & M \ar[r]^x\ar[d] &  M\ar[d]\ar[r]^x & M\ar[r]\ar[d] & \hdots \\
		& \vdots & \vdots & \vdots &
	}
	\]
	We obtain two spectral sequences from this, induced by taking either the horizontal filtration or vertical filtration.  Let us consider the spectral sequence obtained from the horizontal filtration. The first page is given by
	\[
	\xymatrix{
		\hdots \ar[r] & M_y \ar[r]&  M_y\ar[r] & M_y\ar[r] & \hdots \\
		\hdots \ar[r] &  M_y \ar[r] & M_y\ar[r] &  M_y\ar[r] & \hdots \\
		\hdots \ar[r] & M_y \ar[r] &  M_y\ar[r] & M_y\ar[r] & \hdots \\
	}
	\]
	where operator is the induced action of $x$, which we write as $d_1$.  Notice that $h$ acts on $M_y$ because it normalizes $y$, and we have $[h,d_1]=d_1$. The next page is given by
	\[
	\xymatrix{
		\hdots & (M_y)_x &  (M_y)_x & (M_y)_x & \hdots \\
		\hdots \ar[rru] &  (M_y)_x \ar[rru] & (M_y)_x\ar[rru] &  (M_y)_x & \hdots \\
		\hdots \ar[rru] & (M_y)_x \ar[rru] &  (M_y)_x\ar[rru] & (M_y)_x & \hdots \\
	}
	\]
	We call this odd operator $d_{3}$ on $(M_y)_x$, and again $h$ will act on $(M_{y})_{x}$, such that $d_3$ increases the weight of $h$ by $3$.  \
	
	In general on the $n$th page we obtain the same super vector space at every position, which we write as $DS_{y,x}^nM$, and which admits an action of $h$ such that the differential $d_{2n-1}$ increases the $h$-weight by $2n-1$.

	\subsection{The functors \texorpdfstring{$DS_{y,x}^{n}$}{DSxyn} and \texorpdfstring{$DS_{x,y}^{n}$}{DSxyn}}\label{ssec:DS_infty_functors}
	
	By picking out one position in the above spectral sequence, we obtain a sequence of spaces which admit a semisimple action of the even operator $h$, and the action of a square-zero, odd differential $d_{2n-1}$ such that 
	\[
	[h,d_{2n-1}]=(2n-1)d_{2n-1}.
	\]
	Similarly, if we take the spectral sequence associated to the vertical filtration, we obtain a sequence of spaces which admit a semisimple action of $h$ and the action of the square-zero, odd operator $d_{1-2n}$ such that $[h,d_{1-2n}]=(1-2n)d_{1-2n}$.
	
	In general, define functors
	\[   DS_{y,x}^n:\Rep_{\g_{\ol{0}}}\g\l(1|1)\to\Rep\C\langle h,d_{2n-1}\rangle,
	\]
	where $DS_{y,x}^nM$ is the $(0,0)$-position on the $n$th page of the spectral sequence via the horizontal filtration obtained from $M^c$.  Similarly, we have
	\[   DS_{x,y}^n:\Rep_{\g_{\ol{0}}}\g\l(1|1)\to\Rep\C\langle h,d_{1-2n}\rangle,
	\]
	obtained as the $(0,0)$-position on the $n$th page of the spectral sequence via the vertical filtration on $M^c$.
	
	In particular:
	\[
	DS_{y,x}^0M=M^c, \ \ \ \ DS_{y,x}^1M=M_y, \ \ \ \ DS^2_{y,x}M=(M_y)_x.
	\]

	We say these spectral sequences stabilize if the differential vanishes after some page.  In this case we write
	\[
	DS^{\infty}_{y,x}M, \ \ \ DS_{x,y}^\infty M
	\]
	for what they stabilize to, and these spaces will admit a semisimple action of $h$. The spectral sequences clearly stabilize whenever the weights of $M$ under the action of $h$ are bounded; in particular this holds in the case when $M$ is finite-dimensional.  Thus we obtain functors
	\[
	DS_{x,y}^\infty,DS_{y,x}^\infty:\Rep_{\g_{\ol{0}}}\g\l(1|1)\to\Rep\C\langle h\rangle
	\]
	
	The proof of the following theorem, along with all technical details surrounding the spectral sequence, are given in the appendix.
	\begin{thm}
		The functors:
		\begin{itemize}
			\item $DS_{y,x}^n:\Rep_{\g_{\ol{0}}}\g\l(1|1)\to \Rep\C\langle h,d_{2n-1}\rangle$,
			\item $DS_{x,y}^n:\Rep_{\g_{\ol{0}}}\g\l(1|1)\to \Rep\C\langle h,d_{1-2n}\rangle$,
			\item $DS_{x,y}^\infty,DS_{y,x}^\infty:\Rep_{\g_{\ol{0}}}\g\l(1|1)\to \Rep\C\langle h\rangle$,
		\end{itemize}
		are symmetric monoidal functors for all $n\in\N$.
	\end{thm}
	
	\begin{remark}[Caution] Unlike $DS_{y,x}^1=DS_{y}$, which is middle exact, $DS_{y,x}^k$ will not be middle exact for $k>1$.
	\end{remark}
	
	\subsection{Relationship to contragredient duality}\label{ssec:spectral_seq_and_contragred}  Recall the order 4 automorphism $\sigma_{\g\l(1|1)}$ of $\g\l(1|1)$ given by $X \mapsto -X^{st}$, where $X^{st}$ denotes the supertranspose. In particular, $\sigma_{\g\l(1|1)}(x)=\pm y$. 
	
	Write $\sigma_n:\langle h,d_{2n-1}\rangle\to\langle h,d_{1-2n}\rangle$ for the isomorphism of Lie superalgebras given by $h\mapsto -h$ and $d_{2n-1}\mapsto d_{1-2n}$.  Then we may twist a $\C\langle h,d_{2n-1}\rangle$-module by $\sigma_n$ to obtain a $\C\langle h,d_{1-2n}\rangle$-module.
	
	\begin{lemma}\label{lemma twisting DS}
		We have a natural isomorphism of functors
		\[
		DS^i_{x,y}\circ(-)^{\sigma_{\g\l(1|1)}}\cong (-)^{\sigma_n}\circ DS^i_{y,x}.
		\]
	\end{lemma} 
	\begin{proof}
		See Section \ref{sec_appendix_contra_duality} of the appendix.
	\end{proof}

	Now define the contragredient duality functors  $(-)^\vee$ as the composition of dualizing with twisting by either $\sigma_{\g\l(1|1)}$ or $\sigma_h$.  Note \InnaA{that the functor} $(-)^\vee$ has already been defined and discussed for $\Rep\g\l(1|1)$ in \cref{ssec:contra duality}.  Even though for any object $V$ in $\Rep\g_m$ we have an isomorphism $V\cong V^\vee$, this isomorphism is not natural, so it is important not to identify these modules.
	
	\begin{thm}\label{thm commutes contragredient}
		For $n\in\N\cup\{\infty\}$, we have a canonical isomorphism in $\Rep\g_m$:
		\[
		DS^n_{x,y}(M^\vee)\cong (DS^n_{y,x}M)^\vee.
		\]
	\end{thm} 
	\begin{proof}
		Combine the isomorphisms from \cref{prop DS dual} and \cref{lemma twisting DS}.
	\end{proof}	
	
	\subsection{Action of \texorpdfstring{$d_n$}{dn} on indecomposable \texorpdfstring{$\p\g\l(1|1)$}{pgl(1|1)}-modules}\label{ssec:spectral_seq_and_indecomposables}  We now describe the action of $DS_{y,x}^n$ and $DS_{x,y}^n$ on indecomposable $\p\g\l(1|1)$-modules.  For $r\in\C$, we write $\C_r$ for the one-dimensional even $\C\langle h\rangle$-module on which $h$ acts by $r$.
	
	\begin{lemma}\label{action_on_fds}
		\mbox{}
		\begin{enumerate}
			\item $DS_{y,x}^nQ=DS_{x,y}^nQ=0$ for all $n>0$ and for any projective module $Q \in Rep(\p\g\l(1|1))$.
			\item $DS_{y,x}^nW(m)=\C_{-m}$ for all $n>0$ and all $m\in\Z$.
			\item $DS_{x,y}^nW(m)=\C_{m}$ for all $n>0$ for all $m\in\Z$.
			\item $DS_{y,x}^nY(m)=0$ for all $n>0$.
			\item $DS_{x,y}^nX(m)=0$ for all $n>0$.
			\item $DS_{y,x}^nX(m)=\C_{-m+1/2}\oplus \Pi\C_{m-1/2}$ for $n\leq m$, $DS_{y,x}^nX(m)=0$ for $n>m$.  If we write $\C_{-m+1/2}=\langle u_{-m+1/2}\rangle$ and $\Pi\C_{m-1/2}=\langle u_{m-1/2}\rangle$, then $d_i=0$ for $i\neq 0,m$ and $d_m(u_{-m+1/2})=u_{m-1/2}$.
			\item $DS_{x,y}^nY(m)=\C_{-m+1/2}\oplus\Pi\C_{m-1/2}$ for $n\leq m$, $DS_{x,y}^nX(m)=0$ for $n>m$.  If we write $\C_{-m+1/2}=\langle u_{-m+1/2}\rangle$ and $\Pi\C_{m-1/2}=\langle u_{m-1/2}\rangle$, then $d_i=0$ for $i\neq 0,m$ and $d_m(u_{m-1/2})=u_{-m+1/2}$.
		\end{enumerate}
	\end{lemma}
	\begin{proof}
		See Section \ref{sec_appendix_proof_fds}.
	\end{proof}

	\subsection{Tensor product rules for \texorpdfstring{$\Rep \g\l(1|1)$}{Rep(gl(1|1))}}\label{ssec:tens_prod_gl_1_1}
	Using \cref{action_on_fds} and the fact that $DS_{y,x}^r$ and $DS_{x,y}^r$ are symmetric monoidal functors, one can easily determine the tensor products of modules in $\Rep\g\l(1|1)$ up to projective summands (see \cite{GQS07} for the precise formulas which include the projective summands).
	
	In the next proposition, we write $M \cong N \oplus Proj$ whenever $M$ is isomorphic to a direct sum of $N$ with some projective $\g\l(1|1)$-module.
	
	\begin{prop}\label{prop tensor prod gl(1|1)}
		We have the following tensor product rules on $\Rep\g\l(1|1)$ up to a projective summand: 
		\begin{enumerate}
			\item For $m,n\in\Z$ we have $W(m)\otimes W(n)\equiv W(m+n)  \oplus  Proj$;
			\item For $n\in\Z$, $m\geq0$, we have $W(n)\otimes X(m)\cong X(m)_{-n}  \oplus  Proj$;
			\item For $n\in\Z$, $m\geq0$, we have $W(n)\otimes Y(m)\cong Y(m)_{n} \oplus  Proj$;
			\item For $0<m\leq n$, we have $X(m)\otimes X(n)\cong X(m)_{-n+1/2}\oplus\Pi X(m)_{n-1/2} \oplus  Proj$;
			\item For $0<m\leq n$, we have $Y(m)\otimes Y(n)\cong Y(m)_{-n+1/2}\oplus\Pi Y(m)_{n-1/2} \oplus  Proj$;
			\item For $m,n\in\Z$, $X(m)\otimes Y(n)$ is projective.
		\end{enumerate}
	\end{prop}
	\begin{proof}
		In each case we apply $DS_{x,y}^\infty$ and $DS_{y,x}^{\infty}$, using that these are symmetric monoidal functors, and determine the decomposition.
	\end{proof}
	From part (1) of \cref{prop tensor prod gl(1|1)}, one can compute the semisimplification of $\Rep_{\g_{\ol{0}}}\g\l(1|1)$, however we prove this directly later on.
	
	\begin{remark}
		One could determine the full decomposition of the tensor products in \cref{prop tensor prod gl(1|1)} using the $\C\langle h \rangle$-character of the tensor product and the description of the indecomposable projectives of $\p\g\l(1|1)$ given in \cref{ssec:indecomps}.
	\end{remark}
	
	\subsection{Application: tensor products in \texorpdfstring{$\Rep_{\g_{\ol{0}}}\g\l(1|n)$}{Rep(gl(1|n))}}\label{ssec:tens_prod_gl_1_n} 
	
	Let $\g=\g\l(1|n)$. Choose an isotropic root $\alpha$ and a corresponding root subalgebra $\g\l(1|1)\sub\g\l(1|n)$. We write $x\in\g_{\alpha}$, $y\in\g_{-\alpha}$, $c=[x,y]$, and $h$ as usual.  We are interested in the stable category of $\Rep_{\g_{\ol{0}}}\g\l(1|n)$: \InnaA{this is the category obtained from $\Rep_{\g_{\ol{0}}}\g\l(1|n)$ after quotienting by the ideal of morphisms which factor through a projective $\g\l(1|n)$-module}. The blocks of this category, up to parity, are indexed by irreducible representations of $\g\l(n-1)$, see \cite{GS10}.  Write this set as $\II_{n-1}$, \InnaB{and denote by $\BB_V'$ the block corresponding to $V\in\II_{n-1}$. The block $\BB_V'$ is equivalent to the principal block of $\Rep(\g\l(1|1))$ (see \cite{GS10})}.  
	
	Thus for each $V\in\II_{n-1}$, we may consider a sum of blocks $\BB_{V}:=\BB_V'\oplus\Pi\BB_V'$ of $\Rep_{\g_{\ol{0}}}\g\l(1|n)$.  We may index the simple modules of $\BB_{V}$ up to parity with integers, and we write $L_V(n)$ with $n\in\Z$ for the unique simple module in $\BB_V$ such that $DS_xL_V(n)=V$ and the $h$ action on $DS_xL_V(n)$ is multiplication by the scalar $n$.  
	
	Given $n_1,n_2\in\Z$, we set $W_V(n_1;n_2)$ to be the module in $\BB_{V}$ corresponding to the $\g\l(1|1)$-module $W(n_1)_{n_2}$.  Define $X_V(n_1;n_2)$, $Y_V(n_1;n_2)$ when $n_1\in\N$, $n_2\in\Z$ to be the modules corresponding to $X(n_1)_{n_1-1/2+n_2}$ and $Y(n_1)_{n_1-1/2+n_2}$. Because the functors $DS_{x,y}^m$ and $DS_{y,x}^m$ commute with translation functors, we are able to use similar methods as in for instance \cite{GH20}, to obtain the following: 
	\begin{enumerate}
		\item $DS_{x,y}^mW_V(n_1;n_2)=DS_xW_V(n_1;n_2)=V_{n_1+n_2}$ for all $m>0$;
		\item $DS_{y,x}^mW_V(n_1;n_2)=DS_yW_V(n_1;n_2)=V_{n_1-n_2}$ for all $m>0$;
		\item $DS_{x,y}^mX_V(n_1;n_2)=0=DS_{y,x}^mY_V(n_1;n_2)$ for $m>0$;
		\item $DS_{y,x}^mX_V(n_1;n_2)=V_{n_2}\oplus\Pi V_{2n_1-1+n_2}$ for $m\leq n_1$;
		\item $DS_{y,x}^mX_V(n_1;n_2)=V_{n_2}\oplus\Pi V_{2n_1-1+n_2}$ for $m>n_1$;
		\item $DS_{x,y}^mY_V(n_1;n_2)=V_{n_2}\oplus\Pi V_{2n_1-1+n_2}$ for $m\leq n_1$;
		\item $DS_{x,y}^mY_V(n_1;n_2)=V_{n_2}\oplus\Pi V_{2n_1-1+n_2}$ for $m>n_1$.
	\end{enumerate}
	
	For \InnaA{irreducible representations} $V(\lambda),V(\mu), V(\gamma)$ \InnaA{of $\mathfrak{gl}(n-1)$ we denote by $c_{\lambda\mu}^{\gamma}$ the multiplicity of $V(\gamma)$ in $V(\lambda)\otimes V(\mu)$ (the Littlewood-Richardson coefficient): in other words,} $V(\lambda)\otimes V(\mu)=\bigoplus\limits_{\gamma}V(\gamma)^{c_{\lambda\mu}^{\gamma}}$.  
	
	\begin{thm}
		We have the following tensor product relations in the stable category of $\Rep_{\g_{\ol{0}}}\g\l(1|n)$:
		\begin{enumerate}
			\item $W_{V(\lambda)}(n_1;n_2)\otimes W_{V(\mu)}(m_1;m_2)\equiv \bigoplus\limits_{\gamma}W_{V(\gamma)}(n_1+m_1;n_2+m_2)^{c_{\lambda\mu}^{\gamma}}$;
			\item $W_{V(\lambda)}(n_1;n_2)\otimes X_{V(\mu)}(m_1;m_2)\equiv\bigoplus\limits_{\gamma}X_{V(\gamma)}(m_1;m_2-n_1+n_2)^{c_{\lambda\mu}^{\gamma}}$;
			\item $W_{V(\lambda)}(n_1;n_2)\otimes Y_{V(\mu)}(m_1;m_2)\equiv\bigoplus\limits_{\gamma}Y_{V(\gamma)}(m_1;m_2+n_1+n_2)^{c_{\lambda\mu}^{\gamma}}$;
			\item $0<n_1\leq m_1$: 
			\[
			X_{V(\lambda)}(n_1;n_2)\otimes X_{V(\mu)}(m_1;m_2)\equiv\bigoplus\limits_{\gamma}X_{V(\gamma)}(n_1;n_2+m_2)^{c_{\lambda\mu}^{\gamma}}\oplus \Pi X_{V(\gamma)}(n_1;n_2+m_2+2m_1-1)^{c_{\lambda\mu}^{\gamma}};
			\]
			\item $0<n_1\leq m_1$: 
			\[
			Y_{V(\lambda)}(n_1;n_2)\otimes Y_{V(\mu)}(m_1;m_2)\equiv\bigoplus\limits_{\gamma}Y_{V(\gamma)}(n_1;n_2+m_2)^{c_{\lambda\mu}^{\gamma}}\oplus\Pi Y_{V(\gamma)}(n_1;n_2+m_2+2m_1-1)^{c_{\lambda\mu}^{\gamma}}.
			\]
		\end{enumerate}
	\end{thm}
	
	From the above relations, we see that the semisimplification of the category $\Rep_{\g_{\ol{0}}}\g\l(1|n)$ is $\Rep\g_m\times\Rep GL(n-1)\times\Rep(G_m)$. This computation was originally done in \cite{H19}.  
	
	\section{Realizing the semisimplification}\label{sec:realizing_the_ss}
	\subsection{}
	
	Recall the category $\Rep\g_m$ of semisimple representations of $\C\langle h \rangle$. So far we have constructed, for each $n\in\N\cup\{\infty\}$, symmetric monoidal functors 
	\[ DS_{y,x}^n:\Rep_{\g_{\ol{0}}}\g\l(1|1)\to\Rep\C\langle h,d_{2n-1}\rangle,
	\]
	\[ DS_{x,y}^n:\Rep_{\g_{\ol{0}}}\g\l(1|1)\to\Rep\C\langle h,d_{1-2n}\rangle,
	\]
	where $d_{\pm\infty}=0$, and which are interchanged by contragredient duality.  
	For the following, recall that for a $\g\l(1|1)$-module $V$ and $r\in\C$, we write $V_r$ for the tensor product of $V$ by the one dimensional module of character $r\mathbf{str}$.
	\begin{lemma}\label{DS_infty_on_maps}
		Let $m, n \in \Z$, $r,s \in \C$.
		
		\begin{enumerate}
			\item  The $\C\langle h \rangle$-module
			$DS_{x, y}^{\infty}W(n)_r$ is one-dimensional, with $h$ acting by the scalar $n+r$. We have
			\[
			DS_{x, y}^{\infty}\Hom_{\g\l(1|1)}(W(n)_r,W(m)_s)=
			\begin{cases}
				\C &\text{ if }s-r=\InnaA{n-m \geq 0} \\
				0 &\text{ else}
			\end{cases}.
			\]
			\item The $\C\langle h \rangle$-module
			$DS_{y,x}^{\infty}W(n)_r$ is one-dimensional, with $h$ acting by the scalar $-n+r$. We have
			\[
			DS_{y, x}^{\infty}\Hom_{\g\l(1|1)}(W(n)_r,W(m)_s)=	\begin{cases}
				\C  &\text{ if }s-r=\InnaA{m-n\leq 0} \\
				0 &\text{ else}
			\end{cases}. 
			\]
		\end{enumerate}
	\end{lemma}
	
	\begin{proof}
		\InnaA{The statements about the modules $DS_{x,y}^{\infty}W(n)_r$, 	$DS_{y,x}^{\infty}W(n)_r$ follow from \cref{action_on_fds}.  To compute the images of the Hom-spaces, recall from \cref{prop tensor prod gl(1|1)} that
			\begin{align*}
				\Hom_{\g\l(1|1)}(W(n)_r,W(m)_s) &\cong \Hom_{\g\l(1|1)}(W(n),W(m)_{s-r})\cong \Hom_{\g\l(1|1)}(\C, W(-n)\otimes W(m)_s)\\&
				\cong\Hom_{\g\l(1|1)}(\C,W(m-n)_s \oplus Q)
			\end{align*}
			for some projective $\g\l(1|1)$-module $Q$. Since the functors $DS_{x,y}^{\infty}$, $DS_{y,x}^{\infty}$ are symmetric monoidal and send $Q$ to zero, we conclude that it is enough to prove the statements in the case $n=r=0$. In this case the non-zero $\g\l(1|1)$-maps $W(0)\cong \C \to W(m)_s$ are embeddings, so the required statements follow immediately from the definitions of the functors.}
	\end{proof}

	The following diagram illustrates the morphisms between such objects $DS_{y,x}^{\infty}W(n)_r$, $ DS_{y,x}^{\infty}W(n)_r$ on which $h$ acts with integral eigenvalues. We write $\ol{(-)}$ for the images of the indecomposable $\g\l(1|1)$-modules in $\Rep\g_m$; the red (respectively, blue) arrows show the images of $\g\l(1|1)$-morphisms under the functor $DS_{y,x}^{\infty}$ (respectively, $ DS_{y,x}^{\infty}$).
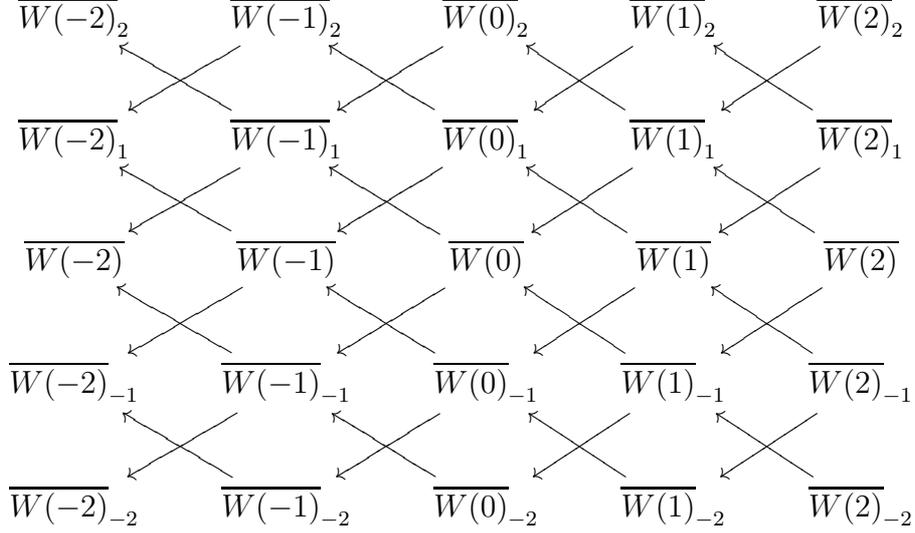
	\[
	\xymatrix{
		\ol{W(-2)}_{2} & \ol{W(-1)}_2\ar@[blue][ld]  & \ol{W(0)}_2 \ar@[blue][dl]& \ol{W(1)}_2 \ar@[blue][dl]& \ol{W(2)}_2\ar@[blue][ld]	\\
		\ol{W(-2)}_{1} & \ol{W(-1)}_{1} \ar@[blue][dl] \ar@[red][lu]& \ol{W(0)}_{1}\ar@[blue][ld] \ar@[red][lu] & \ol{W(1)}_1\ar@[blue][dl] \ar@[red][lu]& \ol{W(2)}_1\ar@[blue][ld] \ar@[red][lu]\\
		\ol{W(-2)} & \ol{W(-1)}\ar@[red][lu]\ar@[blue][dl] & \ol{W(0)}\ar@[red][lu]\ar@[blue][dl] & \ol{W(1)}\ar@[red][lu]\ar@[blue][dl] & \ol{W(2)}\ar@[red][lu]\ar@[blue][dl]\\
		\ol{W(-2)}_{-1} & \ol{W(-1)}_{-1}\ar@[red][lu]\ar@[blue][dl] & \ol{W(0)}_{-1}\ar@[red][lu]\ar@[blue][dl] & \ol{W(1)}_{-1}\ar@[red][lu]\ar@[blue][dl] & \ol{W(2)}_{-1}\ar@[red][lu]\ar@[blue][dl]\\
		\ol{W(-2)}_{-2} & \ol{W(-1)}_{-2}\ar@[red][lu] & \ol{W(0)}_{-2}\ar@[red][lu] & \ol{W(1)}_{-2}\ar@[red][lu]& \ol{W(2)}_{-2}\ar@[red][lu] .
	}
	\]
	
	\begin{remark}
		\InnaA{In particular, although the objects $DS_{x,y}^{\infty}W(n)_r$, $DS_{x,y}^{\infty}W(n+1)_{r-1}$ are isomorphic as $\C\langle h \rangle$-modules, the inverse map is not in the image of $DS_{x,y}^\infty$. Thus we see that the functors $DS_{y,x}^{\infty}$ and $ DS_{y,x}^{\infty}$ are not full.}
	\end{remark}
	
	\subsection{Filtrations on \texorpdfstring{$DS_{x,y}^{\infty}$}{DSxy} and \texorpdfstring{$DS_{y,x}^{\infty}$}{DSyx}} \label{ssec:filtr}
	\subsubsection{Definition}
	\InnaA{The functors $DS_{x,y}^{\infty} V$, $DS_{y,x}^{\infty} V$ as defined before had a significant downside: they sent non-isomorphic indecomposable $\g\l(1|1)$-modules to isomorphic one-dimensional $\C\langle h\rangle$-modules. In order to remedy this, we will now show that these functors are naturally equipped with additional structure. }
	
	We define natural filtrations on our two functors, both of which we denote by $F_{\bullet}$.  Infomrally, the filtration works as follows: for $V\in \Rep_{\g_{\ol{0}}}\g\l(1|1)$, we filter the object $DS_{x,y}^{\infty} V$ using the above diagram by images of indecomposable summands $\ol{W(k)}_{\bullet}$ lying to the left of some vertical line, i.e. such that $k\leq n$ for some corresponding $n$.

	\begin{example}\label{ex:filtr_image_W_n}
		If $V = W(n)_s$, we see from the above diagram that the 
		$(1|0)$-dimensional module $DS_{x,y}^{\infty} W(n)_s$ has a natural filtration, 
		which is essentially "$2$-step": $F_{k} DS_{x,y}^{\infty} W(n)_s = 0$ for $k<n$, 
		and $F_{k} DS_{x,y}^{\infty} W(n)_s = DS_{x,y}^{\infty}W(n)_s$ for $k\geq n$. 
	\end{example}

	More precisely, for $V\in\Rep_{\g_{\ol{0}}}\g\l(1|1)$ and $n\in\Z$, set
	\[
	F_nDS_{x,y}^\infty V:=\im(\,\Hom_{\s\l(1|1)}(W(n),V)\otimes DS_{x,y}^\infty W(n) \to DS_{x,y}^\infty V\,),
	\]
	and
	\[
	F_nDS_{y,x}^\infty V:=\im(\,\Hom_{\s\l(1|1)}(W(n),V)\otimes DS_{y,x}^\infty W(n) \to DS_{y,x}^\infty V\,),
	\]
	
	In other words, $F_nDS_{x,y}^\infty V$ contains all vectors in $DS_{x,y}^\infty V$ which lie in the images of morphisms of the form $DS_{x,y}^\infty f: DS_{x,y}^\infty W(n) \to DS_{x,y}^\infty V$ for some $f\in \Hom_{\s\l(1|1)}(W(n),V)$, and similarly for the second functor.
	
	\begin{remark}
		\InnaA{To avoid dealing with the twists of the $W(n)$'s, which would make the definition more cumbersome, we used maps of $\s\l(1|1)$-module and instead of $\g\l(1|1)$.}
	\end{remark}	
	
	\InnaA{Let us show that this is indeed a filtration:
		\begin{lemma}\label{lem:filtr_well_def}
			For any $V \in \Rep_{\g_{\ol{0}}}\g\l(1|1)$, we have: $F_nDS_{x,y}^\infty V \subset F_{n+1}DS_{x,y}^\infty V$ for any $n$, and \[ \bigcap\limits_n F_n DS_{x,y}^\infty V=0, \ \ \ \ \bigcup\limits_n F_nDS_{x,y}^\infty V =DS_{x,y}^\infty V.\]
			An analogous statement holds for the filtration $F_{\bullet} DS_{y,x}^\infty V $.
		\end{lemma}
		\begin{proof}
			We prove the statements for the functor $DS_{x,y}^\infty$, for the other functor the proof is analogous. 
			
			Fix maps $f_n: W(n+1)_{-1} \to W(n)$ such that  $DS_{x,y}^\infty(f_n) \neq 0$ (so $DS_{x,y}^\infty(f_n)$ is an isomorphism in $\Rep\g_m$). Let $\phi: W(n) \to V$. Then \[\im(DS_{x,y}^\infty(\phi) \circ DS_{x,y}^\infty(f_n)) =\im( DS_{x,y}^\infty(\phi)) \] so $\im DS_{x,y}^\infty(\phi) \subset  F_{n+1}DS_{x,y}^\infty V$, proving the first statement.
			
			Next, decomposing $V$ into a direct sum of indecomposable $\g\l(1|1)$-modules, we see that $F_n DS_{x,y}^\infty(V) \neq 0$ iff $V$ contains a summand of the form $W(k)_r$, $k\leq n$, $r\in \C$. This implies that $\bigcap\limits_n F_n DS_{x,y}^\infty V=0$. Moreover, since $ DS_{x,y}^\infty V$ is a direct sum of $DS_{x,y}^\infty W(k)_r$ where $W(k)_r$ are direct summands of $V$, we conclude that $\bigcup\limits_n F_nDS_{x,y}^\infty V $ is the entire space $DS_{x,y}^\infty V$.
		\end{proof}
	}
	
	\subsubsection{A new target category}
	It is easy to check that if $\varphi:V\to W$ is a morphism of $\g\l(1|1)$-modules, then the map $DS_{x,y}^{\infty}\varphi$ (respectively, $DS_{y,x}^\infty\varphi$) respects the filtrations defined above. 
	
	\InnaA{This implies that the filtrations define functors $\Rep_{\g_{\ol{0}}}\g\l(1|1)\to\Rep\g_m^{fil}$, where $\Rep\g_m^{fil}$ is the category of filtered, semisimple $\C\langle h\rangle$-modules. By abuse of notation, we will denote these functors again by $DS_{x,y}^\infty$ and $DS_{y,x}^\infty$.}

	\begin{lemma}\label{lem:DS_infty_filt_is_SM}
		The functors $DS_{x,y}^\infty,DS_{y,x}^\infty: \Rep_{\g_{\ol{0}}}\g\l(1|1)\to\Rep\g_m^{fil}$ are symmetric monoidal functors.
	\end{lemma}
	\begin{proof}
		We prove the statement for the functor $DS_{x,y}^\infty$, the proof for $DS_{y,x}^\infty$ being analogous. 
		
		To show that $DS_{x,y}^\infty$ is a symmetric monoidal functor, we only need to establish a natural transformation $$DS_{x,y}^\infty V\otimes DS_{x,y}^\infty W \longrightarrow DS_{x,y}^\infty(V\otimes W)$$ and show that it is an isomorphism (the fact that this functor respects the unit and the symmetry morphisms is obvious). Indeed, given $V, V' \in \Rep_{\g_{\ol{0}}}\g\l(1|1)$, the filtered $\C\langle h\rangle$-module $DS_{x,y}^\infty(V) \otimes DS_{x,y}^\infty(V')$ has $n$-th filtration given by
		\[F_n\left(DS_{x,y}^\infty(V) \otimes DS_{x,y}^\infty(V')\right) = \sum_{n_1+n_2= n} F_{n_1} DS_{x,y}^\infty(V) \otimes F_{n_2} DS_{x,y}^\infty(V')\]
		This is a subspace of the vector superspace $DS_{x,y}^\infty(V) \otimes DS_{x,y}^\infty(V')$ which is isomorphic (as a superspace) to $DS_{x,y}^\infty(V \otimes V')$. This subspace is spanned by the subspaces $\im DS_{x,y}^\infty(f) \otimes\im DS_{x,y}^\infty(f')$ for all $f \in \Hom_{\s\l(1|1)}(W(n_1), V)$ and $f' \in \Hom_{\s\l(1|1)}(W(n_2), V')$, $n_1+n_2=n$. Such a pair $(f, f')$ defines a map $f\otimes f': W(n_1)\otimes W(n_2)\to V \otimes V'$. By \cref{prop tensor prod gl(1|1)}, $W(n_1)\otimes W(n_2)$ is isomorphic to a direct sum of $W(n)$ with some projective $\g\l(1|1)$ module. Thus $f\otimes f'$ induces a map $DS_{x,y}^\infty(f\otimes f'): DS_{x,y}^\infty(W(n)) \to DS_{x,y}^\infty(V\otimes V')$. Hence we obtain a natural inclusion 
		\begin{equation}\label{eq:DS_filt_SM}
			F_n\left(DS_{x,y}^\infty(V) \otimes DS_{x,y}^\infty(V') \right) \, \subset \, F_n DS_{x,y}^\infty(V\otimes V')
		\end{equation} making the functor $DS_{x,y}^\infty$ lax monoidal. To show that it is, in fact, strongly monoidal, we need to show that the above inclusion is actually an isomorphism, a statement which it is enough to verify whenever $V, V'$ are indecomposable $\g\l(1|1)$-modules. The only case of interest here is when $V\cong W(k)_r$, $V'\cong W(m)_s$ for some $k, m\in \Z$, $r, s\in \C$. In that case, \cref{ex:filtr_image_W_n} shows that the dimensions of both sides of the inclusion \eqref{eq:DS_filt_SM} are equal, which implies that this inclusion is an isomorphism.
	\end{proof}
	\subsubsection{The associated graded}
	We may now define
	\begin{align}\label{eq:def_ass_gr}
		DS_{x,y}^\infty V[n]&:=F_nDS_{x,y}^\infty V/F_{n-1}DS_{x,y}^\infty V,\\
		DS_{y,x}^\infty V[n]&:=F_nDS_{y,x}^\infty V/F_{n-1}DS_{y,x}^\infty V.
	\end{align}
	
	Define corresponding cofiltrations $F^\bullet$ on $DS_{x,y}^\infty$ and $DS_{y,x}^\infty$ as
	\[
	F^nDS_{x,y}^\infty V=DS_{x,y}^\infty V/F_{-n\InnaA{-1}}DS_{x,y}^\infty V,
	\]
	and
	\[
	F^nDS_{y,x}^\infty V=DS_{y,x}^\infty V/F_{-n\InnaA{-1}}DS_{y,x}^\infty V.
	\]
	Clearly we have canonical isomorphisms:
	\begin{align}\label{eq:alternative_ass_gr}
		\ker(F^{n}DS_{x,y}^\infty V\to F^{n-1}DS_{x,y}^\infty V)&\cong DS_{x,y}^\infty V[-n]\\
		\ker(F^{n}DS_{y,x}^\infty V\to F^{n-1}DS_{y,x}^\infty V)&\cong DS_{y,x}^\infty V[-n]
	\end{align}
	Now we can prove:
	\begin{prop}\label{prop iso contra}
		The natural isomorphism $ DS_{x,y}^\infty V^\vee\to (DS^\infty_{y,x}V)^\vee$ of \cref{thm commutes contragredient} takes the filtration $F_\bullet DS_{x,y}^\infty V^\vee$ to the cofiltration $F^\bullet DS^\infty_{y,x}V$.  In particular it induces a natural isomorphism
		\[
		DS_{x,y}^\infty V^\vee[n]\cong (DS_{y,x}^\infty V[-n])^\vee.
		\]
		\InnaA{where $(-)^{\vee}$ on the left hand side stands for contragredient duality of $\g\l(1|1)$ and on the right hand side stands for contragredient duality of $\C\langle h \rangle$-modules.}
	\end{prop}

	\begin{proof}
		The isomorphism $DS_{x,y}^\infty V^\vee\to (DS^\infty_{y,x}V)^\vee$ gives rise to a natural filtration on $(DS^\infty_{x,y}V)^\vee$, which thus gives a natural cofiltration on $DS^\infty_{y,x}V$.  In order to show this cofiltration agrees with the one defined above, it suffices to check this is so on indecomposables, and argue by naturality.  
		
		Thus consider the indecomposable module $W(n)_r$.  Since $W(n)_r^\vee=W(-n)_{-r}$, we have
		\[
		F_kDS_{x,y}^\infty W(n)_r^\vee= \begin{cases}
			0, & \text{if } k<-n \\
			DS_{x,y}^\infty W(n)_r^\vee, & \text{if }k\geq -n.
		\end{cases}
		\]
		It follows that the induced cofiltration on $DS_{y,x}^\infty W(n)_{r}$ is given by
		\[
		F^kDS_{x,y}^\infty W(n)_r= \begin{cases}
			0, & \text{if } k<-n \\
			DS_{x,y}^\infty W(n)_r, & \text{if }k\geq -n.
		\end{cases}
		\]
		This is exactly the cofiltration defined above.
	\end{proof}

	The following result should be viewed as a concrete realization of the uniqueness of semisimplification up to isomorphism of functors; see \cref{ssec:semisimplification} for more precise statements. 
	
	\begin{prop}\label{prop iso DS_x DS_y}
		For each $n\in\Z$ we have a natural isomorphism of $\C\langle 
		h\rangle$-modules
		\[
		DS^\infty_{x,y}V[n]\to DS^\infty_{y,x}V_{2n}[n]
		\]
		Here on the right hand side we have the $(2n)$-twist of $V$ by the Berezinian.
	\end{prop}
	
	\begin{proof}
		\InnaB{In this proof all the isomorphisms will be up to parity shift.}
		Let $v\in DS_{x,y}^\infty V[n]$. There exists a summand $W\cong W(n)$ of $V$ such that $v$ lifts to the lowest (respectively, the highest) weight vector of $W$ if $n\leq0$ (respectively, $n>0$).  We denote this vector by $\tilde{v}$.  Now, within this summand $W$ there is a natural isomorphism between the highest and lower weight spaces; namely, the maps $x,y$ define isomorphisms between the weight spaces, so we may repeatedly uniquely lift and project. Explicitly, our isomorphism is given by:
		\[
		\tilde{v}\mapsto (xy^{-1})^n\tilde{v} \text{ for } n\leq 0, \ \ \ \ \tilde{v}\mapsto(x^{-1}y)^n\tilde{v} \text{ for } \ n\geq0.
		\]
		In this way we may map $\tilde{v}$ to the highest (respectively, the lowest) weight vector, which we call $w$. Then $w$ defines an element in $DS_{y,x}^\infty V[n]$ with weight shifted by $-2n$; thus to make this map $h$-equivariant, we need to twist by $2n$. Once we show the above procedure is well-defined, it is clear that it is natural (i.e. respects $\C\langle h\rangle$-maps), and so we will be done.
		
		To check that it is well-defined, we need to check what happens if we choose another summand $W'\cong W(n)$ of $V$ in which $v$ lifts to the lowest (respectively highest) weight vector $\tilde{v}'$.  We will explain what happens in the case of $n=-m\leq0$, with the case $n>0$ following from applying \cref{prop iso contra}.  Notice that since $DS_{x,y}^\infty W(n)=DS_xW(n)$, our choice of $\tilde{v},\tilde{v}'$ have that $\tilde{v}-\tilde{v}'\in\operatorname{im}x$.
		
		Under our setup, $\tilde{v},\tilde{v}'$ lie in the socles \InnaA{of $W$, $W'$ respectively}, and give rise to well-defined bases
		\[
		w_{-m}=\tilde{v},w_{1-m},\dots,w_{m}, \  \text{ and } \ w_{-m}'=\tilde{v}',w_{1-m}',\dots,w_{m}'
		\]
		of $W$ and $W'$ respectively. These bases are obtained by using the fact that $x$ and $y$ define isomorphisms between weight spaces. Here these bases have the property that $w_{-m+2i}$ lies in the socle of $W$ for all $i$, $yw_{-m+2i+1}=w_{-m+2i}$, $xw_{-m+2i+1}=w_{-m+2i+2}$ for all $0\leq i\leq m-1$, and similarly for the corresponding vectors in $W'$. 
		
		We want to show that $w_{m}$ and $w_{m}'$ admit equivalent projections 
		to $DS_{y,x}^\infty V[n]_{2n}$; assume for a contradiction that this is 
		not the case.  
		
		Take $z_i=w_i-w_i'$; then we claim that $z_m$ will be the highest weight 
		vector of a submodule of $V$ isomorphic to either $X(j)_k$  for some 
		$j>0$ and $k\in \Z$ or to $W(j)_k$ for some $j<n$, $k\in\Z$. Indeed, 
		\InnaA{$xz_{m-1}=z_m \neq 0$ by our assumption, so} we must have $z_{m-1}\neq0$; 
		if $z_{m-2}=yz_{m-1}=0$ then $z_{m-1},z_m$ define a submodule isomorphic to 
		$X(1)_k$; and if $z_{m-2}\neq0$, then we may continue the argument. 
		So that either for some $i$ we have \InnaA{$z_{m-2i}=0$ and then} 
		$z_{m-2i+1},\dots,z_m$ span a submodule isomorphic to $X(i)_k$ for 
		some $k \in \C$, or \InnaA{$z_{m-2i}\neq 0$ for all $i=0, 1 \ldots, 2m$. In the 
			latter case,} $z_{-m},\dots,z_m$ span a submodule isomorphic to $W(n)_k$ for 
		some $k$. But as noted above, $z_{-m}=\tilde{v}-\tilde{v}'\in\im x$, so we must have that 
		$z_{-m}=xz_{-m-1}$ for some element $z_{-m-1}\in V$; if $yz_{-m-1}=0$, then we 
		obtain $X(m+1)_k$ for some $k$, and if $yz_{-m-1}\neq0$ then we obtain a 
		submodule isomorphic to $W(-m-1)_k$.
		
		Let $M$ denote this indecomposable submodule of $V$ containing $z_m$ 
		that is isomorphic to either $X(j)_k$ for some $k$ or $W(j)_k$ with $j<n$ for 
		some $k$.  Since by assumption $z_{m}$ 
		defines a nontrivial element in $DS_{y,x}^\infty V[n]_{2n}$, there must 
		exist a split copy of $W(n)$ in $V$ such that $z_{m}$ is the highest weight 
		vector.  Thus we obtain a nontrivial map $M\to W(n)$ which is nonzero on 
		$z_{m}$; but by \cref{lemma 
			morphism} such a map must be zero on the socle of $M$, which contains $z_m$, 
		which leads to a contradiction. This proves the required statement.
		
	\end{proof}
	
	We now obtain the main theorem. 
	\begin{thm}\label{theorem culmative}
		We have a natural isomorphism of $\C\langle h\rangle$-modules, respecting tensor products:
		\[
		DS_{x,y}^{\infty}V^\vee[n]\cong \left( DS_{x,y}^\infty V_{-2n}[-n]\right)^\vee.
		\]
		\InnaA{Here on the right hand side we have the $(-2n)$-twist of $V$ by the Berezinian, and $(-)^{\vee}$ on the right hand side again stands for contragredient duality of $\C\langle h \rangle$-modules.}
	\end{thm}
	
	\subsection{Categorical viewpoint}\label{ssec:semisimplification}  
	\subsubsection{}
	\InnaA{Recall that the filtrations defined functors $DS_{x,y}^\infty,DS_{y,x}^\infty:\Rep_{\g_{\ol{0}}}\g\l(1|1)\to\Rep\g_m^{fil}$ into the category $\Rep\g_m^{fil}$ of filtered semisimple $\C\langle h\rangle$-modules. Unlike the functors we previously defined into $\Rep\g_m$, these functors turn out to be very nice:
		\begin{thm}\label{thrm equiv fil}
			The functors $DS_{x,y}^\infty,DS_{y,x}^\infty:\Rep_{\g_{\ol{0}}}\g\l(1|1)\to\Rep\g_m^{fil}$ are essentially surjective, full, symmetric monoidal functors.
		\end{thm}
		\begin{proof}
			\InnaA{We prove the theorem for the functor $DS_{x,y}^\infty$, the proof for $DS_{y,x}^\infty$ being analogous. We have already seen in \cref{lem:DS_infty_filt_is_SM} that $DS_{x,y}^\infty$ is symmetric monoidal, so we only need to show that it is full and essentially surjective.

				To prove that $DS_{x,y}^\infty$ is essentially surjective, it is enough to show that any indecomposable filtered semisimple $\C\langle h\rangle$-module is obtained as the image of some $W(n)_r$ or its parity shift. Indeed, any such indecomposable is one-dimensional, and thus isomorphic (up to parity shift) to some  filtered $\C\langle h\rangle$-module $ DS_{x,y}^\infty W(n)_r$ as in \cref{ex:filtr_image_W_n}. }
			
			We now show that $DS_{x,y}^\infty$ is full. It is enough to check that the map
			\[ DS_{x,y}^\infty: \Hom_{\g\l(1|1)}(V, V') \longrightarrow \Hom_{\Rep\g_m^{fil}}(DS_{x,y}^\infty V, DS_{x,y}^\infty V') \] is surjective when $V, V'$ are indecomposable $\g\l(1|1)$-modules, and more specifically, both of the form $W(n)_r$ for some $n, r$ (otherwise $DS_{x,y}^\infty V =0$ for indecomposable $V$). In that case, by \cref{DS_infty_on_maps}, the map
			\[ DS_{x,y}^\infty: \Hom_{\g\l(1|1)}(W(n)_r, W(m)_s) \longrightarrow \Hom_{\Rep\g_m^{fil}}(DS_{x,y}^\infty W(n)_r, DS_{x,y}^\infty W(m)_s) \] has a one-dimensional image if $s-r = n-m \geq 0$, and $0$ otherwise. On the other hand, $DS_{x,y}^\infty W(n)_r$ is a filtered $\C\langle h\rangle$-module of weight $n+r$ with filtration \[F_{k}DS_{x,y}^\infty W(n)_r = \begin{cases}
				0 &\text{ if } k<n\\
				DS_{x,y}^\infty W(n)_r &\text{ if } k \geq n\\
			\end{cases}\] so the space $\Hom_{\Rep\g_m^{fil}}(DS_{x,y}^\infty W(n)_r, DS_{x,y}^\infty W(m)_s)$ is $0$ if $n<m$ and one-dimensional otherwise. This proves that $DS_{x,y}^\infty$ is full.
	\end{proof}}

	\subsubsection{Semisimplifications} Recall from \cite{AKO02} the notion of the semisimplification of a rigid symmetric monoidal category. Namely, let $\AA$ be a rigid symmetric monoidal category, and let $\NN$ be the monoidal ideal of $\AA$ given by all the negligible morphisms: 
	\InnaA{\[\NN = \{f: X \to Y \in Mor(\AA)\;: \;\; \forall g: Y\to X, \, Tr(g\circ f)=0\}.\]} Then one may consider the quotient $S: \AA \to \AA^{ss}:=\AA/\NN$. We call the pair $(S, \AA^{ss})$ the semisimplification of $\AA$; \InnaA{if $\AA$ satisfies the property that the trace of any nilpotent endomorphism is zero, then $\AA^{ss}$ is a semisimple category.  }
	
	By \cite[Proposition 2.3.4]{AKO02}, if $\SS$ is any semisimple rigid symmetric monoidal category for which there exists a full, essentially surjective symmetric monoidal functor $S':\AA\to\SS$, then there exists an isomorphism of categories $F:\AA^{ss}\to\SS$ such that the two functors \InnaA{$S, S'\circ F:\AA\to\SS$} are isomorphic.
	
	Consider the category $\Gr\Rep\g_m:= \Rep\g_m\times\Rep\G_m$ and the functor $\Gr:\Rep\g_m^{fil}\to \Gr\Rep\g_m$ sending a filtered $\C\langle h\rangle$-module to its associated graded. \InnaA{This functor is clearly symmetric monoidal, as well as full and essentially surjective, making the pair $(\Gr, \Gr\Rep\g_m)$ a semisimplification of the rigid symmetric monoidal category $\Rep\g_m^{fil}$.}
	
	We define functors $DS_{x,y}^{ss},DS_{x,y}^{ss}:\Rep_{\g_{\ol{0}}}\g\l(1|1)\to\Gr\Rep\g_m$ by
	\[
	DS_{x,y}^{ss}:=\Gr\circ DS_{x,y}^{\infty}, \ \ \ \ DS_{y,x}^{ss}:=\Gr\circ DS_{y,x}^{\infty}.
	\] 
	The grading can be written explicitly as follows:
	\[
	DS_{x,y}^{ss}V=\bigoplus\limits_{n\in\Z}DS_{x,y}^\infty V[n], \ \ \ DS_{y,x}^{ss}V=\bigoplus\limits_{n\in\Z}DS_{y,x}^\infty V[n],
	\]
	where $V[n]$ lies in $\Rep\g_m$.   
	
	The following is an immediate consequence from \cref{thrm equiv fil}.
	\begin{cor}\label{prop functors ss}
		The functors $ DS_{x,y}^{ss},DS_{y,x}^{ss}$ are full, essentially surjective symmetric monoidal functors. \InnaA{This makes the pairs $(DS_{x,y}^{ss}, \Rep\g_m\times\Rep\G_m)$, $(DS_{y,x}^{ss}, \Rep\g_m\times\Rep\G_m)$ into (isomorphic) semisimplifications of the category $\Rep_{\g_{\ol{0}}}\g\l(1|1)$.}
	\end{cor}
	
	\begin{remark}
		It was first proven in \cite{H19} that the semisimplification of $\Rep GL(1|1)$ is $\Rep\G_m\times\G_m$.
	\end{remark}
	
	Proposition \ref{prop functors ss} along with \cite[Proposition 2.3.4 ]{AKO02} give an immediate (albeit less satisfying) proof of \cref{prop iso DS_x DS_y}. Stated in terms of our semisimplification functors, the result becomes:
	
	\begin{cor}
		Consider the the automorphism $\phi$ of $\C\langle h\rangle\times\C$ satisfying $\phi(h,z)=(h+2z,z)$ where $\C$ represents the Lie algebra of $\G_m$. This defines an autofunctor $\Phi_{\phi}$ on $\Rep\g_m\times \Rep\G_m$ given by twisting by $\phi$. Then we have a natural isomorphism of symmetric monoidal functors:
		\[
		DS_{x,y}^{ss}\cong \Phi_{\phi}\circ DS_{y,x}^{ss}.
		\]
	\end{cor}
	
	\InnaA{This isomorphism is realized explicitly by the natural isomorphisms given in \cref{prop iso DS_x DS_y}.}
	
	\subsubsection{Summarizing diagram}
	We may summarize the situation with the following diagram:
	
	
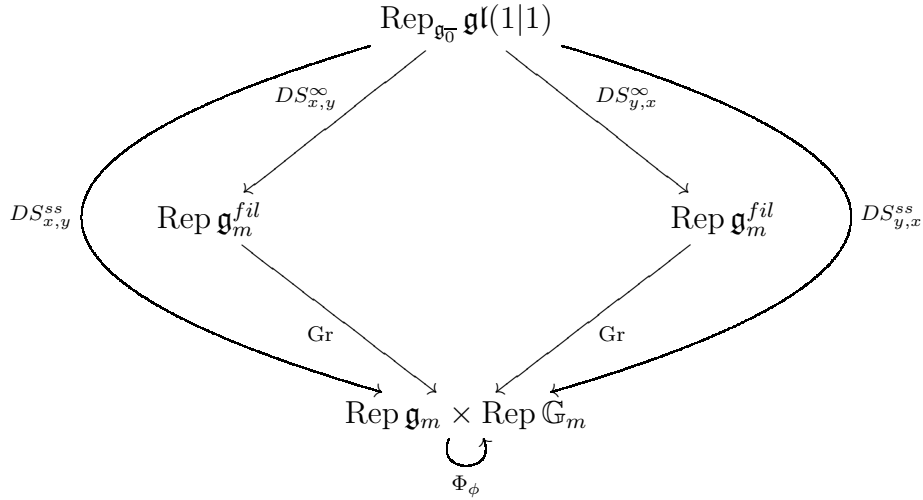
\[
	\xymatrix{ & \Rep_{\g_{\ol{0}}}\g\l(1|1)  \ar@/^12pc/[dddd]^{DS_{y,x}^{ss}} \ar@/_12pc/[dddd]_{DS_{x,y}^{ss}} \ar[ddl]_{DS_{x,y}^{\infty}} \ar[ddr]^{DS_{y,x}^\infty} & \\ 
		&& \\
		\Rep\g_m^{fil}\ar[ddr]_{\operatorname{Gr}} & & \Rep\g_m^{fil} \ar[ddl]^{\operatorname{Gr}}\\
		&&\\
		& \Rep\g_m\times\Rep\G_m  \ar@(dl,dr)_{\Phi_{\phi}} &
	}
	\]
	The left and right triangles in this diagram commute, and the square gives an isomorphism of symmetric monoidal functors after composing with $\Phi_{\phi}$.

	\subsection{Extending the story to the functor \texorpdfstring{$DS_{x+y}$}{DS{x+y}}}\label{ssec:DSx+y}
	
	\InnaA{We now consider the functor $DS_{x+y}: \Rep_{\g_{\ol{0}}}(\g\l(1|1)) \to \SVec$.}
	\begin{lemma}\label{lem:ds_x_plus_y_indec}
		If $M$ is an indecomposable $\g\l(1|1)$-module, then $DS_{x+y}M\neq 0$ if and only if $M\cong W(n)_r$ for some $n\in\Z$ and $r\in\C$.  Furthermore, $\dim DS_{x+y}W(n)_r = (1|0)$.
	\end{lemma}
	\begin{proof}
		For projective $M$ and for $M =W(n)_r$ this follows from \cref{lemma DS on indecomps}. For $M$ of the form $X(n)_r$ or $Y(n)_r$, this is a straightforward computation.
	\end{proof}
	
	\begin{lemma}\label{lem:ds_x_plus_y_homs}
		Let $n, m\in \Z$, $r,s \in \C$. 
		
		\InnaA{If $n\geq m$ and $s-r \in \{m-n,m-n+2,\dots,n-m\}$} then
		\[DS_{x+y}\Hom_{\g\l(1|1)}(W(n)_r,W(m)_s)= \C.
		\]
		Otherwise $	DS_{x+y}\Hom_{\g\l(1|1)}(W(n)_r,W(m)_s) = 0$.
	\end{lemma}
	\begin{proof}
		Similarly to the proof of \cref{DS_infty_on_maps}, we may reduce the problem to the question of computing the spaces $DS_{x+y}\Hom_{\g\l(1|1)}(W(0),W(m)_s) \neq 0$ for all $m\in \Z, s\in\C$. Since $\dim DS_{x+y}W(m)_s=(1|0)$, we know that this space is of dimension at most $1$.
		From here one can compute directly.
	\end{proof}
	We give the local picture of the morphisms given by $DS_{x+y}(f)$ where $f$ is a map between $\g\l(1|1)$-modules of the form $W(n)_r$ with $n, r\in \Z$. We write $\ol{W(n)_r}:=DS_{x+y}W(n)_r$. 
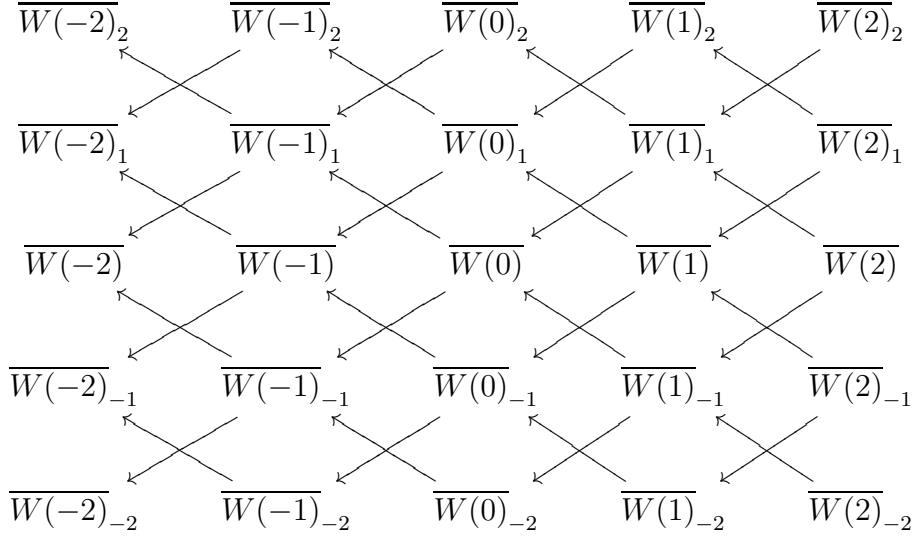
	\[
	\xymatrix{
		\ol{W(-2)}_{2} & \ol{W(-1)}_2\ar[ld]  & \ol{W(0)}_2 \ar[dl]& \ol{W(1)}_2 \ar[dl]& \ol{W(2)}_2\ar[ld]	\\
		\ol{W(-2)}_{1} & \ol{W(-1)}_{1} \ar[dl] \ar[lu]& \ol{W(0)}_{1}\ar[ld] \ar[lu] & \ol{W(1)}_1\ar[dl] \ar[lu]& \ol{W(2)}_1\ar[ld] \ar[lu]\\
		\ol{W(-2)} & \ol{W(-1)}\ar[lu]\ar[dl] & \ol{W(0)}\ar[lu]\ar[dl] & \ol{W(1)}\ar[lu]\ar[dl] & \ol{W(2)}\ar[lu]\ar[dl]\\
		\ol{W(-2)}_{-1} & \ol{W(-1)}_{-1}\ar[lu]\ar[dl] & \ol{W(0)}_{-1}\ar[lu]\ar[dl] & \ol{W(1)}_{-1}\ar[lu]\ar[dl] & \ol{W(2)}_{-1}\ar[lu]\ar[dl]\\
		\ol{W(-2)}_{-2} & \ol{W(-1)}_{-2}\ar[lu] & \ol{W(0)}_{-2}\ar[lu] & \ol{W(1)}_{-2}\ar[lu]& \ol{W(2)}_{-2}\ar[lu] .
	}
	\]
	
	\begin{remark}
		One observes that up in the above picture there are two \InnaA{``blocks of morphisms''}: the block containing $DS_{x+y}W(0)$, and the block containing $DS_{x+y}W(0)_{1}$.  This is reflected in the action of \InnaA{the group $\Z/2\Z$ on the functor} $DS_{x+y}$, which acts by $1$ on the ``block'' containing $DS_{x+y}W(0)$ and by $-1$ on the ``block'' containing $DS_{x+y}W(0)_1$.
	\end{remark}

	Let \InnaA{$\mathbf{a_1}=(1,-1)$, $\mathbf{a_2}=(1,1)$} be elements of 
	$\Z\times\C$.  Define an order on $\Z\times\C$ by setting 
	\InnaA{$\mathbf{b_1}\leq \mathbf{b_2}$} if 
	$\mathbf{b_2}-\mathbf{b_1}=j\mathbf{a_1}+k\mathbf{a_2}$ for some $j,k\in\N$. 
	Next, let $\operatorname{Fil}^{\Z\times\C}$ be the category of pairs 
	$(V, \{V^{\mathbf b}\}_{b\in \Z\times \C}$) such that 
	\begin{itemize}
		\item $V$ is a finite-dimensional super vector space; 
		\item for every $\mathbf{b}\in\Z\times\C$, $V^{\mathbf{b}}\sub V$ is a subspace;
		\item if $\mathbf{b_1}, \mathbf{b_2}\in\Z\times\C$ such that $\mathbf{b_1}\leq \mathbf{b_2}$, then $V^{\mathbf{b_1}}\sub V^{\mathbf{b_2}}$;
		\item $\cup_{\mathbf{b}\in\Z\times\C} V^{\mathbf{b}} = V$.
	\end{itemize}
	
	\InnaA{The morphisms in $\operatorname{Fil}^{\Z\times\C}$ will be morphisms of vector superspaces $\phi:V\to W$ such that $\phi(V^\mathbf{b}) \subset W^\mathbf{b}$ for each $\mathbf{b}\in \Z\times \C$. This category has an obvious symmetric monoidal structure, $(V, \{V^\mathbf{b}\}_{\mathbf{b}}) \otimes (W, \{W^\mathbf{b}\}_{\mathbf{b}}):= (V\otimes W, \{(V\otimes W)^\mathbf{b}\}_{\mathbf{b}})$, where \[\forall\, \mathbf{b}\in \Z\times \C, \;\;\;\;\;\; (V\otimes W)^\mathbf{b} := \sum_{\mathbf{b}'+\mathbf{b}'' = \mathbf{b}} V^{\mathbf{b}'} \otimes W^{\mathbf{b}''} .\]}

	Then $\operatorname{Fil}^{\Z\times\C}$ is a Karoubian monoidal category whose indecomposables are \InnaA{indexed (up to parity shift) by $v\in \Z\times\C$. The indecomposable object corresponding to $v\in \Z\times\C$ is $\C_\mathbf{v}:=(\C, \{\C^{\mathbf{b}}_{\mathbf{v}}\}_{\mathbf{b}})$, where $\C^{\mathbf{b}}_{\textbf{v}}:=\C$ if $\mathbf{b}\geq \mathbf{v}$ and $\C^{\textbf{b}}_{\textbf{v}}=0$ otherwise.} 
	
	Following the same ideas as in \cref{ssec:filtr}, for any $\g\l(1|1)$-module $V$ we may give $\InnaA{\ol{V}:=}DS_{x+y}V$ the structure of an object in $\operatorname{Fil}^{\Z\times\C}$: \InnaA{for any $ \textbf{b} = (n, r)\in \Z\times \C$, \[\ol{V}^\textbf{b} := \im\left( \, \Hom_{\g\l(1|1)} (W(n)_r, V)\otimes DS_{x+y} W(n)_r \to \ol{V}\, \right).\]
		The vectors in the subspace $\ol{V}^\textbf{b} \subset \ol{V}$ are exactly those which lie in the images of morphisms of the form $DS_{x+y}(f): DS_{x+y} W(n)_r \to \ol{V}$ for some $f \in \Hom_{\g\l(1|1)} (W(n)_r, V)$.
		
		The same arguments as in \cref{lem:filtr_well_def} show that this endows $\ol{V}$ with the structure of an object in $\operatorname{Fil}^{\Z\times\C}$. This construction is clearly functorial, defining a functor 
		\[DS_{x+y}: \Rep_{\g_{\ol{0}}} \g\l(1|1)\to \operatorname{Fil}^{\Z\times\C}.\]
		
		\begin{example}\label{ex:ds_x_plus_y_indec}
			Let $V:=W(m)_s$, $m\in \Z, s\in \C$. Set $\mathbf{v}:=(m, s)$. 
			
			Then $\ol{V}:= DS_{x+y}W(n)_r$ is a $(1|0)$-dimensional vector superspace, and by \cref{lem:ds_x_plus_y_homs}, we have: $DS_{x+y}\Hom_{\g\l(1|1)} (W(n)_r, W(m)_s) =\C$ iff $(n, r) \geq (m, s)$ in the above order. Thus we have: 
			$$ \ol{V}^{\textbf{b}} =\C^b_{\mathbf{v}} := \begin{cases}
				\C &\text{ if } \textbf{b}\geq \mathbf{v} \\ 0 &\text{ else }
			\end{cases}$$ for any $\textbf{b}\in \Z\times \C$. In other words, $$\left(DS_{x+y}W(n)_r, \left\{\left(DS_{x+y}W(n)_r\right)^{\mathbf{b}}\right\}_{\mathbf{b}}\right) \cong \C_{\mathbf{v}}.$$ 
		\end{example}
	}
	The following result becomes a natural extension of the ideas used in \cref{thrm equiv fil}.
	
	\begin{prop}
		The functor $DS_{x+y}: \Rep_{\g_{\ol{0}}} \g\l(1|1)\to \operatorname{Fil}^{\Z\times\C}$ is an essentially surjective, full, symmetric monoidal functor.  
	\end{prop}
	\begin{proof}	
		\InnaA{
			The proof that this functor is symmetric monoidal is a direct analogy of the proof of \cref{lem:DS_infty_filt_is_SM}.
			
			From \cref{ex:ds_x_plus_y_indec} we see that for any indecomposable object $\C_{\mathbf{v}} \in \operatorname{Fil}^{\Z\times\C}$ there exists a corresponding indecomposable object in $\Rep_{\g_{\ol{0}}} \g\l(1|1)$ which is sent to $\C_{\mathbf{v}}$. Thus the functor is essentially surjective. 
			
			To see that this functor is full, we only need to check that for any indecomposable $V, W \in \Rep_{\g_{\ol{0}}} \g\l(1|1)$, the map 
			\[DS_{x+y}: \Hom_{\g\l(1|1)}(V, W) \to \Hom_{\operatorname{Fil}^{\Z\times\C}} (DS_{x+y}V, DS_{x+y} W) \] is surjective. By \cref{lem:ds_x_plus_y_indec}, it is enough to show this when $V=W(n)_r$, $W=W(m)_s$ for $n, m\in \Z$, $r,s\in \C$. In that case, by \cref{ex:ds_x_plus_y_indec} we have: $DS_{x+y} W(n)_r \cong \C_{(n, r)}$, $ DS_{x+y} W(m)_s \cong \C_{(m, s)}$, 
			\[\Hom_{\operatorname{Fil}^{\Z\times\C}} (\C_{(n, r)}, \C_{(m, s)}) =\begin{cases}
				\C &\text{ if } (n, r)\geq (m, s) \\ 0 &\text{ else }
			\end{cases}.\] So by \cref{lem:ds_x_plus_y_homs}, \[\Hom_{\operatorname{Fil}^{\Z\times\C}} (\C_{(n, r)}, \C_{(m, s)})  \cong DS_{x+y}  \Hom_{\g\l(1|1)}(W(n)_r, W(m)_s).\] This completes the proof of the proposition.}
	\end{proof}

	\InnaA{Consider the vector superspace}
	\[
	\Gr_{a_i}(V^\bullet)=\bigoplus\limits_{r\in\C} \left.\left(\sum\limits_{k \InnaA{\in \Z}}V^{(0,r)+ka_i}\right)\middle/\left(\sum\limits_{k\InnaA{\in \Z}}V^{(-2,r)+ka_i}\right)\right..
	\]
	\InnaA{To see that this is well defined, recall that $V^{\mathbf{b}} \subset V^{\mathbf{b} + (2,0)}$ for any $\mathbf{b}\in \Z\times \C$, since $(2,0)=\mathbf{a}_1 + \mathbf{a}_2$.}

	For a pictorial explanation \InnaA{of what $Gr_{a_1}(V^\bullet)$ looks like}, we consider the following diagram, where we consider \InnaA{only the direct summand corresponding to} $r=0$: the colored portions (in blue and red) correspond to $\sum\limits_{k \InnaA{\in \Z}}V^{(0,r)+ka_i}$, while the red portion alone corresponds to $\sum\limits_{k\InnaA{\in \Z}}V^{(-2,r)+ka_i}$, which we quotient out by. 
	
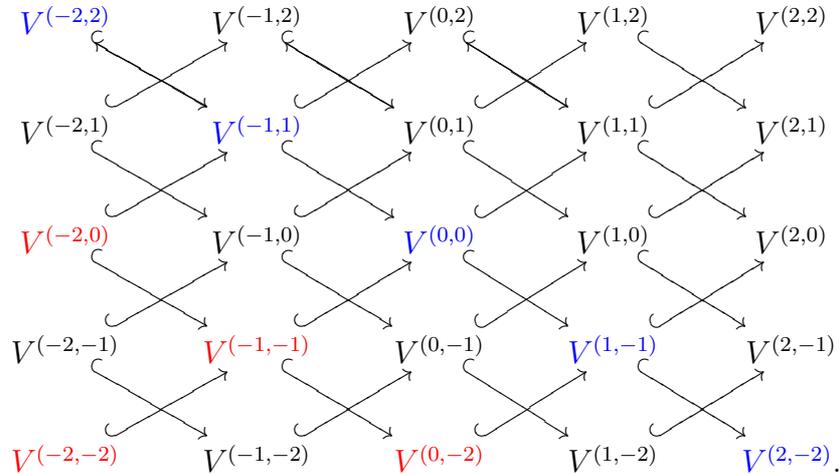
\[
	\xymatrix{
		\color{blue}V^{(-2,2)}\color{black} \ar@{^{(}->}[dr] & V^{(-1,2)}\ar@{^{(}->}[dr]  & V^{(0,2)} \ar@{^{(}->}[dr]& V^{(1,2)} \ar@{^{(}->}[dr]& V^{(2,2)}	\\
		V^{(-2,1)}\ar@{^{(}->}[dr]\ar@{^{(}->}[ur]& \color{blue}V^{(-1,1)}\color{black} \ar@{^{(}->}[ur]\ar@{^{(}->}[dr] \ar[lu]& V^{(0,1)}\ar@{^{(}->}[dr]\ar@{^{(}->}[ur] \ar[lu] & V^{(1,1)}\ar@{^{(}->}[dr]\ar@{^{(}->}[ur] \ar[lu]& V^{(2,1)}\\
		\color{red}V^{(-2,0)}\color{black} \ar@{^{(}->}[dr] \ar@{^{(}->}[ur]& V^{(-1,0)}\ar@{^{(}->}[dr]\ar@{^{(}->}[ur] & \color{blue}V^{(0,0)}\color{black} \ar@{^{(}->}[dr]\ar@{^{(}->}[ur] & V^{(1,0)}\ar@{^{(}->}[dr]\ar@{^{(}->}[ur] & V^{(2,0)}\\
		V^{(-2,-1)} \ar@{^{(}->}[dr]\ar@{^{(}->}[ur] & \color{red}V^{(-1,-1)}\color{black}\ar@{^{(}->}[dr]\ar@{^{(}->}[ur] & V^{(0,-1)}\ar@{^{(}->}[dr]\ar@{^{(}->}[ur] & \color{blue}V^{(1,-1)}\color{black} \ar@{^{(}->}[dr]\ar@{^{(}->}[ur] & V^{(2,-1)}\\
		\color{red}V^{(-2,-2)}\color{black} 
		\ar@{^{(}->}[ur]& V^{(-1,-2)}\ar@{^{(}->}[ur]& \color{red}V^{(0,-2)}\color{black}\ar@{^{(}->}[ur] & V^{(1,-2)}\ar@{^{(}->}[ur]& \color{blue}V^{(2,-2)}\color{black} .
	}
	\]
	
	\InnaA{This allows us to define two functors $\Gr_{a_1},\Gr_{a_2}:\operatorname{Fil}^{\Z\times\C}\to\Rep\g_m^{fil}$, where} the $h$ action \InnaA{on $\Gr_{a_i}(V^\bullet)$} is given by the above \InnaA{$\C$-}grading, and the filtration on each eigenspace of $h$ is given by
	\[
	F_{j}\left.\left[\left.\left(\sum\limits_{k\InnaA{\in \Z}}V^{(0,r)+ka_i}\right)\middle/\left(\sum\limits_{k\InnaA{\in \Z}}V^{(-2,r)+ka_i}\right)\right.\right]\right.=\left.V^{(0,r)\InnaA{+}ja_i}\middle/\left.\left[\left(\sum\limits_{k\in\Z}V^{(-2,r)+ka_i}\right)\cap V^{(0,r)\InnaA{+}ja_i}\right.\right]\right.
	\]
	
	\begin{thm}
		We have natural isomorphisms of symmetric monoidal functors:
		\[
		\Gr_{a_1}\circ DS_{x+y}\simeq DS_{x,y}^{\infty}, \  \ \ \ \Gr_{a_2}\circ DS_{x+y}\simeq DS_{y,x}^{\infty}.
		\]
	\end{thm}
	
	\begin{proof}
		\InnaA{For this proof, whenever we have a semisimple $\C\langle h\rangle$-module $M$, we will write $M_{(a)}$ for the $h$-eigenspace corresponding to the eigenvalue $a$.}
		
		We prove these isomorphisms in the case when $M$ has integral eigenvalues for $h$, with the general case following from twisting by the appropriate multiples of the Berezinian character.  Consider the double complex 
		\[
		\xymatrix{
			& \vdots \ar[d] & \vdots \ar[d] & \vdots \ar[d] & \\
			\hdots \ar[r] & M_{(0)} \ar[r]^x\ar[d]^y &  M_{(1)}\ar[d]^y\ar[r]^x & M_{(2)}\ar[r]\ar[d]^y & \hdots \\
			\hdots \ar[r] &  M_{(-1)} \ar[r]^x\ar[d]^y & M_{(0)}\ar[d]^y\ar[r]^x &  M_{(1)}\ar[r]\ar[d]^y & \hdots \\
			\hdots \ar[r] & M_{(-2)} \ar[r]^x\ar[d] &  M_{(-1)}\ar[d]\ar[r]^x & M_{(0)}\ar[r]\ar[d] & \hdots \\
			& \vdots & \vdots & \vdots &
		}
		\]
		The total complex $C^\bullet$ of this double complex is
		\[
		\cdots\to M^-\xto{x+y} M^+\xto{x+y} M^-\to \cdots.
		\]
		where $\InnaA{C^k=}M^+:=\sum\limits_{j\in2\Z}M_{(j)}$ \InnaA{for any $k\in 
			2\Z$} and $\InnaA{C^k=}M^-:=\sum\limits_{j\in2\Z+1}M_{(j)}$ \InnaA{for any 
			$k\in 2\Z+1$}. In particular, $C^0=M^+$ and $C^1=M^-$.
		
		This double complex induces two spectral sequences in the usual way, via the horizontal or vertical filtrations, or equivalently, via a choice of $x$ or $y$. 
		
		We consider here the vertical filtration, so that we take cohomology in $x$ first, \InnaB{to prove the isomorphism $\Gr_{a_1}\circ DS_{x+y}\simeq DS_{x,y}^{\infty}$. 
			
			The other filtration gives an analogous argument proving the second isomorphism in the statement of the theorem.
			
			Our vertical filtration gives us a spectral sequence whose $n$-th page has $(DS^n_{x,y}M)_{(i+j)}$ at the 
			$(i,j)$-position.} If $M$ is 
		finite-dimensional, then this spectral sequence is regular (see \cite[Theorem 5.5.10]{W}), and thus weakly
		converges. The statement of convergence gives us exactly our natural isomorphism 
		of functors; we spell out what this precisely means for this spectral sequence.
		
		Consider the filtration on our total complex given by $C^\bullet_{\leq 
			n}=\sum\limits_k C^{k}_{\leq (k+n)}$, where for an $h$-module $V$, we set 
		$V_{\leq \ell}=\sum\limits_{i\leq\ell}V_{(i)}$. 
		
		Then the differential preserves this filtration, and since this filtration is 
		complete, it induces a natural filtration $F_kH^i(C^\bullet)$ on the cohomology 
		via taking the image of the filtration in cohomology.
		
		By definition of weak convergence of our spectral sequence we have
		functorial isomorphisms: 
		\[
		F_kH^0(C^\bullet)/F_{k-1}H^0(C^\bullet)\cong (DS_{x,y}^{\infty}M)_{(2k)}, \ \ \ F_kH^1(C^\bullet)/F_{k-1}H^1(C^\bullet)\cong (DS_{x,y}^{\infty}M)_{(2k+1)}
		\]
		This gives the desired natural isomorphisms.  We now need to check that the 
		filtered objects $F_kH^0(C^\bullet)$ and $F_kH^1(C^\bullet)$ agree with the filtration on $DS_{x+y}(M)$ with respect to $\mathbf{a_1}$ that we are taking the associated graded of.  However, by naturality, it suffices to check this on indecomposables alone.  If we consider $M=W(n)_a$ with $a+n$ even, then we obtain 
		$F_kH^0(C^\bullet)=H^0(C^\bullet)$ if $k\geq n+a$, and is 0 otherwise. If 
		$M=W(n)_a$ with $a+n$ odd, then $F_kH^1(C^\bullet)=H^1(C^\bullet)$ if $k\geq 
		n+a+1$, and is 0 otherwise.  This exactly corresponds to the filtration which 
		induces the associated graded $Gr_{a_1}$, as desired.	 
	\end{proof}

	Our diagram from \cref{ssec:semisimplification} can now be enhanced:
	\[
	\xymatrix{ & \Rep_{\g_{\ol{0}}}\g\l(1|1) \ar[dd]^{DS_{x+y}}  \ar@/^12pc/[ddddd]^{DS_{y,x}^{ss}} \ar@/_12pc/[ddddd]_{DS_{x,y}^{ss}} \ar@/_1pc/[dddl]_{DS_{x,y}^{\infty}} \ar@/^1pc/[dddr]^{DS_{y,x}^\infty} & \\ 
		&&\\
		& \operatorname{Fil}^{\Z\times\C}\ar[dr]^{\Gr_{a_2}} \ar[dl]_{\Gr_{a_1}} & \\
		\Rep\g_m^{fil}\ar[ddr]_{\operatorname{Gr}} & & \Rep\g_m^{fil} \ar[ddl]^{\operatorname{Gr}}\\
		&&\\
		& \Rep\g_m\times\Rep\G_m  \ar@(dl,dr)_{\Phi_{\phi}} &
	}
	\]
	
	\section{Applications}\label{sec:applications}
	
	\subsection{Lie superalgebra \texorpdfstring{$(\p)\g\l(1|1)$}{(p)gl(1|1)}-subalgebras}\label{ssec:gl(1|1)}	
	Suppose that $\g$ is an arbitrary Lie superalgebra, and that we have map $\g\l(1|1)\to\g$ which is either an embedding or has kernel spanned by $c$.  Let $\k$ be a subalgebra of $\g$ which commutes with the image of $\g\l(1|1)$ in $\g$.  If $\g$ is a finite type Kac-Moody Lie superalgebra or $\q(n)$, we will be interested in the case when $\g\l(1|1)$ is a diagonal copy of inside a product of root subalgebras isomorphic to $\g\l(1|1)$ (see explanation below), and $\k$ is an embedded copy of $\g_x$.
	
	\begin{lemma} Let $M$ be a finite-dimensional $\g$-module which is semisimple over $\g\l(1|1)_{\ol{0}}$.
		\begin{enumerate}
			\item For each $n\in\N\cup\{\infty\}$, $DS_{x,y}^nM$ and $DS_{y,x}^n M$ naturally have the structure of $\k$-modules such that the differentials $d_n(M)$ define $\k$-equivariant odd differentials.
			\item The spaces $DS_{x,y}^{ss}M$, and $DS_{y,x}^{ss}M$ naturally have the structure of a $\k\times\C\langle h\rangle$-module. 
		\end{enumerate}
	\end{lemma}
	
	\begin{proof}
		This follows from the fact that $[\g\l(1|1),\k]=0$ and thus $\k$ defines $\g\l(1|1)$-equivariant morphisms on $M$.
	\end{proof}
	
	\begin{cor}
		We have natural isomorphisms of $\k\times\C\langle h\rangle$-modules:
		\[
		DS^\infty_{x,y}V[n]\to DS^\infty_{y,x}V[n]_{2n}
		\]
	\end{cor}
	
	\subsection{Finite-type Kac-Moody and queer Lie superalgebras}\label{ssec:KM queer one}  If $\g$ is a finite-type Kac Moody Lie superalgebra or is the queer Lie superalgebra $\q(n)$, then $\g$ admits a Chevalley automorphism $\sigma=\sigma_{\g}$ which acts by $(-1)$ on a maximal even torus and satisfies $\sigma^2=\delta$, where $\delta(u)=(-1)^{\ol{u}}u$ is the grading automorphism on $\g$. 
	
	Now let $\alpha_1,\dots,\alpha_k$ be $k$ linearly independent roots such that $\dim(\g_{\alpha})_{\ol{1}}\neq0$, and such that $\alpha_i\pm\alpha_j$ is not a root for $i\neq j$.  Choose nonzero root vectors $x_i\in\g_{\alpha_i}$ and $y_i\in\g_{-\alpha_i}$, and set $x=\sum\limits_i x_i$, $y=\sum\limits_iy_i$.  Finally let $c=[x,y]$.  Now choose a semisimple element $h\in\g_{\ol{0}}$ such that $[h,x]=x$, $[h,y]=-y$, and $\beta(h)=0$ for all roots $\beta$ such that $\beta\pm\alpha_{i}$ is not a root.  Then we have constructed a subalgebra of $\g$ isomorphic to $\g\l(1|1)$.
	
	\begin{definition}\label{def_gl(11)_root_sub}
		We call a subalgebra isomorphic to $\g\l(1|1)$ as constructed above a diagonal $\g\l(1|1)$ subalgebra.
	\end{definition}
	
	In this setup, by \cite{GHSS22} we may embed $\g_x=\g_y$ into $\g$ such that it commutes with our subalgebra $\g\l(1|1)$.
	
	\begin{cor}\label{cor_X_Y_no_show}
		Let $\g$ be a finite-type Kac-Moody Lie superalgebra or $\q(n)$, and choose a diagonal $\g\l(1|1)$ subalgebra.  Assume that $L$ is a simple finite-dimensional module, and if $\g=\q(n)$ also assume that its highest weight is integral: then the restriction to $\g\l(1|1)$ contains no copies of $X(n)$ or $Y(n)$.
	\end{cor}
	
	See \cref{ssec:half integral} for a reminder on the meaning of integral weights for $\q(n)$.
	
	\begin{proof}		
		We prove that in each case the spectral sequences collapse on the first page, which is clearly equivalent.
		
		Observe that for each spectral sequence, the differential defines an odd, $\g_x$-equivariant endomorphism on each page.  For the Kac-Moody case, by \cite{HW21}, \cite{GH20}, and \cite{GHSS22}, $DS_xL\cong DS_yL$ are pure, meaning that $[DS_xL:L'][DS_xL:\Pi L']=0$ for all simple $\g_x$-modules $L'$.  Thus all differentials must vanish after the 0th page, and the spectral sequence collapses on the 1st page. 
		
		For the $\q(n)$ case we use the computations of \cite{GS22}.  Given any simple $\g_x=\g_y$-module $L'$ of integral weight, it is shown there that for any two composition factors of $DS_xL$ or $DS_yL$ isomorphic to $L'$, the difference of their $h$-weights must be even.  Indeed, this follows from the fact that $DS_{x}L$ is a subquotient of $DS_{x_1}\circ\dots\circ DS_{x_r}(L)$, and we know the statement for the latter module. However, as was noted in \cref{ssec:DS_infty_functors}, the differential $d_n$ in each case has weight $2n-1$ or $1-2n$, so $d_n$ must be $0$ after the zeroth page.
	\end{proof}
	
	\subsection{Half-integral simple modules for \texorpdfstring{$\q(n)$}{q(n)}}\label{ssec:half integral} \subsubsection{} Let $\h\sub\b\sub\q(n)$ be a choice of Cartan subalgebra and Borel subalgebra for $\q(n)$, so that $\t=\h_{\ol{0}}$ is a maximal torus of $\q(n)_{\ol{0}}$.  Then there is a basis $\epsilon_1,\dots,\epsilon_n$ of $\t$ such that the dominant integral weights for $\q(n)$ with respect to this basis are given by:
	\[
	\lambda=(a_1,\dots,a_n), \ \ a_i\in\C \text{ such that } a_i-a_{i+1}\in\N,
	\]
	and if $a_i=a_{i+1}$ then $a_i=0$.  We say a weight $\lambda$ is integral if $a_i\in\Z$ for all $i$, and we say $\lambda$ is half-integral if $a_i+1/2\in\Z$ for all $i$. We say arbitrary module for $\q(n)$ is of half-integral weight if all of its composition factors are half-integral. If $\lambda$ is neither integral nor half-integral then $L(\lambda)$ is projective, so we ignore this case (see for instance \cite{CW12}).
	
	In contrast to integral weight modules for $\q(n)$ as explained in \cref{cor_X_Y_no_show}, half-integral simple modules may contain 'Z'-modules; in fact we have the following:
	\begin{lemma}\label{lemma W no show half int}
		If $V$ is a half-integral weight module for $\q(n)$, and $\g\l(1|1)\sub\q(n)$ is a diagonal subalgebra as in Definition \ref{def_gl(11)_root_sub}, then the restriction of $V$ to $\g\l(1|1)$ contains no $W$-submodules as direct summands.
	\end{lemma}
	\begin{proof}
		Because half-integral modules are projective over the Cartan subalgebra of $\q(n)$, we must have that $DS_{x+y}V=0$.  From this the statement follows.
	\end{proof}
	
	The above lemma tells us that our spectral sequence will have interesting terms for half-integral modules; we now compute them in the case when $V=L$ is a simple module and $\g\l(1|1)$ is a root subalgebra, i.e. we only use one root in our construction of $\g\l(1|1)$.
	
	\subsubsection{Arc diagrams} We recall from \cite{GS22} the arc diagrams associated to half-integral weights.  First of all, the simple half-integral weight modules for $\q(n)$ are indexed by their highest weights, which are given by a strictly decreasing sequence of half-integers, i.e.  $\lambda=(a_1/2,\dots,a_n/2)$ where $a_1/2>a_2/2>\dots>a_n/2$ and $a_1,\dots,a_n$ are odd integers.  Given a highest weight $\lambda=(a_1/2,\dots,a_n/2)$, we associate a weight diagram to it, which will be map $f_{\lambda}:\frac{2\N+1}{2}\to\{\circ,>,<,\times\}$, i.e. a map from elements of the form $a/2$ for $a$ positive odd integer.  It is defined as follows:
	\begin{itemize}
		\item If $2b\neq \pm a_i$ for any $i$, then we declare $f_{\lambda}(b)=\circ$;
		\item if $2b=a_i$ and $2b\neq -a_j$ for any $j$, then $f_{\lambda}(b)=>$;
		\item if $2b=-a_{i}$ for some $i$ and $2b\neq a_j$ for any $j$, then $f_{\lambda}(b)=<$;
		\item finally, if $2b=a_i,-a_j$ for some $i,j$ then $f_{\lambda}(b)=\times$.  
		
	\end{itemize} 
	We visualize $f_{\lambda}$ with its graph; for example, if 
	\[
	\lambda=(15/2,13/2,5/2,1/2,-1/2,-3/2,-5/2,-15/2),
	\]
	then $f_{\lambda}$ looks as follows:
	\[
	\begin{tikzpicture}
		\draw (0,0) -- (5,0);
		\draw (0.5,0) node[label=center:{\large $\times$}] {};
		\draw (1,0) node[label=center:{\large $<$}] {};
		\draw (1.5,0) node[label=center:{\large $\times$}] {};
		\draw (2,0)  circle(3pt);
		\draw (2.5,0)  circle(3pt);
		\draw (3,0)  circle(3pt);
		\draw (3.5,0) node[label=center:{\large $>$}] {};
		\draw (4,0) node[label=center:{\large $\times$}] {};
		\draw (4.5,0)  circle(3pt);
	\end{tikzpicture}
	\]
	The above association defines a bijection between half-integral highest weights for $\q(n)$ and weight diagrams with the $r$ symbols $>$, $s$ symbols $<$, and $t$ symbols $\times$ such that $2t+r+s=n$.
	
	Associated to $f_{\lambda}$ we define an arc diagram, which consists of the above weight diagram along with arcs connecting each symbol $\times$ to a symbol $\circ$ to the right of it, such that (1) the arcs do not intersect and (2) no symbol $\circ$ lies underneath an arc.  This uniquely specifies the arc diagram; for example, associated to the the weight diagram above we obtain the arc diagram:
	
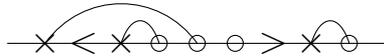
\[
	\begin{tikzpicture}
		\draw (0,0) -- (5,0);
		\draw (0.5,0) node[label=center:{\large $\times$}] {};
		\draw (1,0) node[label=center:{\large $<$}] {};
		\draw (1.5,0) node[label=center:{\large $\times$}] {};
		\draw (2,0)  circle(3pt);
		\draw (1.5,0) .. controls (1.625,0.4) and (1.875,0.4) .. (2,0);
		\draw (2.5,0)  circle(3pt);
		\draw (0.5,0) .. controls (1,0.7) and (2,0.7) .. (2.5,0);
		\draw (3,0)  circle(3pt);
		\draw (3.5,0) node[label=center:{\large $>$}] {};
		\draw (4,0) node[label=center:{\large $\times$}] {};
		\draw (4.5,0)  circle(3pt);
		\draw (4,0) .. controls (4.125,0.4) and (4.375,0.4) .. (4.5,0);
	\end{tikzpicture}
	\]
	Before we can state our theorem, we need some terminology.  We say that an arc is maximal if it does not lie under another arc.  Maximal arcs have the special property that if they are removed along with the symbol $\times$ to which they are attached, an arc diagram for $\q(n-2)$ is obtained.  
	
	We say that a position $\circ$ in an arc diagram is free if it is not the end of any arc.  Given a dominant half-integral weight $\lambda$ and a half-integer $n/2$, define $\ell(\lambda,n/2)$ to be the number of free positions to the left of $n/2$ in the arc diagram of $\lambda$ (see \cite{GS22} for examples).
	
	We caution that in the following theorem we only consider a root subalgebra $\g\l(1|1)$, which is less general than the setting considered in \cref{lemma W no show half int} and \cref{cor_X_Y_no_show}. 
	\begin{thm}
		Let $\lambda$ be a dominant half-integral for $\q(n)$, and $\mu$ a dominant half-integral weight for $\q(n-2)$.  Let $\g\l(1|1)\sub\q(n)$ be a root subalgebra, and let $x,y,h,c$ be its generators.  Then:
		\begin{enumerate}
			\item[(1)] $DS_{x,y}^{k}L(\lambda)$ is semisimple;
			\item[(2)] Either $[DS_{x,y}^{k}L(\lambda):L(\mu)]=0$ or $L(\mu)$ appears in $DS_{x,y}^{k}L(\mu)$ with multiplicity $(1|1)$;
			\item[(3)] We have $[DS_{x,y}^{k}L(\lambda):L(\mu)]\neq0$ if and only if the arc diagram of $\mu$ is obtained from the arc diagram of $\lambda$ by removing a maximal arc such that if the $\times$ end of the arc lies at $j/2$, then we have $k\leq\ell(\lambda,j/2)+1$.
		\end{enumerate}
	\end{thm}
	
	\begin{proof}
		The proof is essentially identical to the computation of $DS_x$ given in \cite{GS22}, and we will explain it using the language from that article.  Namely, because $DS_{x,y}^k$ is a symmetric monoidal functor, and it takes the standard module to the standard module, it will commute with translation functors. Further, it commutes with the operation of removing core symbols for stable weights, because this is given by taking the eigenspace of a semisimple element $z$ which commutes with both $x$ and $y$.  Thus following the algorithm described in \cite{GS22}, we obtain the formula
		\[
		[DS^{k}_{x,y}(L):L(\mu)]=\dim DS_{x,y}^kL(\lambda'),
		\]
		where $\lambda'$ is the weight for $\q_2$ whose arc diagram is obtained by 'shrinking' all arcs of $\mu$ within the arc diagram of $\lambda$.  In this way we reduce the computation to the case of $\q_2$.  Now one may use that the simple module $L(\frac{2n-1}{2}(\epsilon_1-\epsilon_2)))$ decomposes over $\g\l(1|1)$ as $X(n)_{-n/2}\oplus\Pi Y(n)_{-n/2}$.  From this the statement follows.
	\end{proof}
	
	\subsection{Contragredient duality}\label{ssec:applications_contragr}  We continue to let $\g$ denote a finite-type Kac-Moody Lie superalgebra or $\q(n)$.  As previously noted, these Lie superalgebras all admit Chevalley automorphisms $\sigma=\sigma_{\g}$.  Thus we may define contragredient duality functors $V\mapsto V^\vee$ on the category of finite-dimensional modules.  If $V$ is finite-dimensional, we have a canonical isomorphism $V\cong (V^\vee)^\vee$.
	
	Let $\g\l(1|1)\sub\g$ be a subalgebra constructed as in \cref{ssec:KM queer one}.  Let us further assume now that it is stable under $\sigma_{\g}$.   Let $\k$ be a root subalgebra of $\g$ commuting with $\g\l(1|1)$ such that if $\g_{\alpha}\sub\k$ then $\g_{-\alpha}\sub\k$.  In particular this implies that $\k$ is stable under $\sigma_{\g}$.
	
	\begin{thm}\label{theorem contra duality ss}
		We have a natural isomorphism of $\k\times\C\langle h\rangle$-modules:
		\[
		DS_{x,y}^{\infty}V^\vee[n]\cong \left( DS_{x,y}^\infty V[-n]_{-2n}\right)^\vee.
		\]
		where the outer $\vee$ on the right hand side denotes the contragredient duality functor on $\k\times\C\langle h\rangle$-modules which is induced by the restriction of $\sigma$ to this subalgebra.
	\end{thm}
	\begin{proof}
		This follows from the naturality of the isomorphism in \cref{theorem culmative}.
	\end{proof}
	
	\subsection{The category \texorpdfstring{$\Rep_{\g_{\ol{0}}}^+(\g)$}{Rep(g)}}\label{ssec:the_plus_subcategory} Given an abelian symmetric monoidal category $\CC$, let $\CC^+$ denote the Karoubian symmetric monoidal subcategory generated by all simple objects. In other words, the objects of $\CC^+$ are the direct summands of arbitrary tensor products of simple modules. We write $\Rep_{\g_{\ol{0}}}^+(\g):=(\Rep_{\g_{\ol{0}}}(\g))^+$. If $\g=\g\l(1|1)$, then $\Rep_{\g_{\ol{0}}}^+(\g)$ has indecomposables given by all projective modules along with the modules $W(0)_r$ for any $r\in\C$.  These categories (and in particular their semisimplifications) have been studied in the case of $\g=\g\l(m|n)$, see \cite{HW18}.
	\begin{cor}\label{mult free}
		In addition to the hypotheses of \cref{theorem contra duality ss}, assume that $L$ is simple and $DS_xL$ is a multiplicity-free $\k$-module.  Then the restriction of $L$ to $\g\l(1|1)$ lies in $\Rep_{\g_{\ol{0}}}^+(\g\l(1|1))$.
	\end{cor}
	\begin{proof}
		By \cref{cor_X_Y_no_show}, no `Z'-modules appear.  Since $L$ is simple we have $L^\vee\cong L$, and since $DS_{x,y}^\infty L(-n)$ is a subquotient of $DS_xL$ as a $\k$-module, we have $(DS_{x,y}^\infty L(-n))^\vee\cong DS_{x,y}^\infty L(-n)$ as $\k$-modules.  Thus ignoring the $h$-action, \cref{theorem contra duality ss} becomes:
		\[
		DS_{x,y}^{\infty}L(n)\cong DS_{x,y}^\infty L(-n)_{-2n}.
		\]
		Now suppose that for some $n\neq0$ we had $DS_{x,y}^{\infty}L(n)\neq0$, and let $L'$ be a composition factor of it as a $\k$-module.  Then $L'$ must be a composition factor of $DS_{x,y}^\infty L(-n)$.  However then $L$ would need to appear with multiplicity greater than one in $DS_{x,y}^\infty L$; since this is a subquotient of $DS_xL$, we obtain a contradiction.
	\end{proof}
	
	\begin{thm}\label{theorem branching}
		If $L$ is a simple module over $\g\l(m|n)$, $\o\s\p(2m|2n)$, then the restriction of $L$ to a root subalgebra $\g\l(1|1)$ lies in $\Rep_{\g_{\ol{0}}}^+(\g\l(1|1))$.
	\end{thm}
	\begin{proof}
		We use \cref{mult free} in the case of $\k=\g_{x}$ is an embedded subalgebra.
		
		We obtain the cases of $\g\l(m|n)$ and blocks of $\o\s\p(2m|2n)$ which are equivalent to the principal block of $\o\s\p(2k|2k)$ for some $k$, because $DS_xL$ is multiplicity-free in these cases (see \cite{HW21} for the $\g\l$ case and \cite{GH20} for the $\o\s\p$ case).  
		
		The blocks of $\o\s\p(2m|2n)$ that are not equivalent to the principal block of $\o\s\p(2k|2k)$ for any $k$ are instead equivalent to the principal block of $\o\s\p(2k+2|2k)$ for some $k$.	Let $\BB$ be a block equivalent to the principal block of $\o\s\p(2k+2|2k)$, and let $L$ be a simple module.  Then to show that $\operatorname{Res}_{\g\l(1|1)}L\in \Rep_{\g_{\ol{0}}}^+(\g\l(1|1))$, we may apply translation functors to assume that $L$ is stable.  Once doing this, we may find a simple $\o\s\p(2m|2n)$-module $L'$ lying in a block $\BB'$ equivalent to the principal block of $\o\s\p(2k|2k)$, such that $L$ is a direct summand of $T_{\BB'}^{\BB}L'$, where $T_{\BB'}^{\BB}$ is the translation functor from $\BB'$ to $\BB$ (see \cite{GS10}, or Sec. 7.6.3 of \cite{GH20}).  
		
		Now we know already that if we restrict $L'$ to $\g\l(1|1)$ it will lie in $\Rep_{\g_{\ol{0}}}^+(\g\l(1|1))$.  Since the standard module for $\o\s\p(m|2n)$ also has this property, $T_{\BB'}^{\BB}L'$ will inherit this property from $L'$, meaning that $L$ will also have this property, finishing the proof.  
	\end{proof}
	\begin{remark}
		We conjecture that \cref{theorem branching} extends to $\o\s\p(2m+1|2n)$.  However the result does not extend to $\q(n)$; in fact for $\q(2)$ already one see that if $\alpha$ is the unique simple positive root, then for $n\in\N$ we have
		\[
		\operatorname{Res}_{\g\l(1|1)}^{Q(2)}L(n\alpha)=W(n)\oplus\Pi W(-n).
		\]
	\end{remark}
	
	\section{Appendix: the spectral sequences}	
	
	Here we give all details behind the construction of our functors obtained via the spectral sequences.  Because we are not \emph{exactly} using a full spectral sequence, rather only one position of it, we give all details below.  Further, we will need the notation to give a clear proof of Lemma \ref{action_on_fds}.
	
	\subsection{Explicit terms of the sequence}\label{ssec:spectral_seq_def}  Since every position on the spectral sequence introduced in Section \ref{sec:spectral_seq} at any page is the same, we fix the position $(0,0)$ and study what is happening there.  We give an explicit description of the spectral sequence at this position with inspiration from Ravi Vakil's definition in terms of $(p,q)$-strips, see \cite{V08}. We will use this to obtain formulas for boundary maps, which will be odd endomorphisms on each page.
	
	\begin{remark}
		The spectral sequences we construct do not require $M$ to be finite-dimensional; however in order to have convergence this will be a necessary condition.
	\end{remark}  
	
	For each $r\in\Z_{\geq0}$ we will define spaces $B_r(M)=B_r$, $Z_r(M)=Z_r$ and $E_r(M)=E_r$.  The terms of the spectral sequence will be given by $E_r$, while $Z_r$ and $B_r$ will denote respectively the cycles and boundaries in $E_{r-1}$ (we set $E_{-1}:=M$ with $d_{-1}=0$).
	
	\begin{definition}\label{def:r_chain}
		\InnaA{Let $r\geq 1$. A sequence $(v_1,\dots,v_r)$, where $v_1,\dots,v_r \in M$, is called an {\it $r$-chain} in $M$ if it satisfies:  
			\[yv_i=xv_{i+1}\text{ for all } i\leq r-1, \, yv_r=0.\]}
	\end{definition}
	\InnaA{Thus an $r$-chain is a sequence of the form
		\[
		\xymatrix@C=10pt{&v_r\ar@{|->}[dr]^x \ar@{|->}[dl]_y &  &{v_{r-1}} \ar@{|->}[dr]^x \ar@{|->}[dl]_y & & \hdots &{v_2} \ar@{|->}[dr]^x & &v_{1}  \ar@{|->}[dl]_y \\
			{\quad \quad 0\quad \quad} &	&{xv_r=yv_{r-1}} & &{xv_{r-1}=yv_{r-2}} &\hdots  & &{xv_2=yv_1}  }
		\]}
	Now set $Z_0=M$, and for $r>0$ set:
	\begin{equation}\label{eqdef:Z_r}
		Z_r=\{v_r \in M\,:\,\exists \text{ an } r\text{-chain } (v_1,\dots,v_{r-1}, v_r)\}.
	\end{equation}

	In particular we see that $Z_r\sub\ker y$ for $r>0$, and $Z_1=\ker y$.  Now set $B_0=0$, $B_1=\im y$ and for $r>1$ define
	\begin{equation}\label{eqdef:B_r}
		B_r=\im y+\{xw_1\,:\,\exists \text{ an } (r-1)\text{-chain } (w_1,\dots,w_{r-1}) \}.
	\end{equation}
	\InnaA{
		\begin{lemma}\label{lem:prop_Z_r}
			We have: 
			\begin{enumerate}
				\item $Ker(x)\cap Ker(y) \subset Z_r$ for all $r\geq 0$.
				\item $Z_{r+1} \subset Z_r$ for all $r\geq 0$.
				\item $B_r \subset B_{r+1}$ for all $r\geq 0$.
				\item $B_{r'}\sub Z_r$ for all $r, r'\geq 0$.
			\end{enumerate}
		\end{lemma}
		\begin{proof}
			For the first statement, if $v \in Ker(x)\cap Ker(y)$, then the sequence $(0, \ldots, 0, v_r:=v)$ satisfies the conditions in Definition \eqref{eqdef:Z_r} so $v\in Z_r$.
			
			For the second statement, let $v_{r+1} \in Z_{r+1}$ and let $(v_1, v_2, \ldots, v_r, v_{r+1})$ be an $(r+1)$-chain. Then $(v_2, \ldots, v_r, v_{r+1})$ is an $r$-chain so $v_{r+1}\in Z_{r+1}$.
			
			Next, to show that $B_r \subset B_{r+1}$, we only need to check that for $xw_1 \in B_r$ as in Definition \eqref{eqdef:B_r}, we also have $xw_1 \in B_{r+1}$. Indeed, the $(r-1)$-chain $(w_1, \ldots, w_{r-1})$  can be extended to an $r$-chain $(w_1, \ldots, w_{r-1}, 0)$ and thus $xw_1 \in B_{r+1}$ as well.
			
			For the last statement, if we take an element $yv\in\im y$ then we have $y(yv)=0$ and $x(yv)=-y(xv)$ so we may set $v_{r-1}=-xv$. Then $xv_{r-1}=-x^2v=0$, so we can take $v_{i}=0$ for $i<r-1$. Then $(v_1, \ldots, v_{r-1}, v_r:=v)$ is an $r$-chain so $v\in Z_r$.
			
			On the other hand, if we have an $(r-1)$-chain $(w_1, \ldots, w_{r-1})$ then $yw_1=xw_2$, so that $x(xw_1)=0$, and $yxw_1=-xyv_1=-x^2w_2=0$. Thus $xw_1 \in Ker(x) \cap Ker(y)$ and so $xw_1 \in Z_r$ for any $r$. 
		\end{proof}
		
		Thus we have two chains of subspaces, one increasing and one decreasing:
		\[ 0= B_0 \subset B_1 \subset B_2 \subset \ldots \subset \ldots \subset Z_2 \subset Z_1 \subset Z_0 = M.\]
	}
	
	We define
	\[
	E_r(M)=E_r:=Z_r/B_r.
	\]
	Given $v\in Z_r$, we write $\ol{v}$ for its projection to $E_r$.  The following is clear.  
	\begin{lemma}
		The definitions of $Z_r,B_r$, and $E_r$ are functorial in $M$.
	\end{lemma}
	\begin{example}
		\InnaA{	$E_0 = M$, $E_1 = M_y$, $E_2 = DS_x (DS_y M)$.}
	\end{example}
	\subsection{The differential}\label{ssec:spectral_seq_differential} 
	\InnaA{
		We now define a differential $d_r(M)=d_r$ on $E_r$. First, we define a \InnaA{parity-shifting} map $\widetilde{d}_r:Z_r \to E_r$ by setting $\widetilde{d}_0:=y$, and for $r>0$ and any $v_r\in Z_r$ we let
		\[
		\widetilde{d}_r(\ol{v_r}):=\ol{xv_1}.
		\] where $(v_1, \ldots, v_r)$ is an $r$-chain. }
	
	\InnaA{	First, we show that this map is well defined:
		\begin{lemma}
			We have: $xv_1 \in Z_r$, and $\ol{xv_1}$ does not depend on the choice of sequence $(v_1, \ldots, v_r)$.		
		\end{lemma}
		\begin{proof}
			First we check that $xv_1\in Z_r$. Indeed, $xv_1 \in B_{r+1}$ and by Lemma \ref{lem:prop_Z_r} we have: $B_{r+1} \subset Z_r$. So $xv_1\in Z_r$. 
			
			Secondly, suppose that we have $2$ $r$-chains $(v_1,\dots,v_{r-1}, v_r)$ and $(v_1',\dots,v_{r-1}', v_r)$. Let $w_i=v_i-v_i'$.  Then we have that $yw_{r-1}=0$, $xw_{r-1}=yw_{r-2},\dots,xw_2=yw_1$ so $(w_1, \ldots, v_{w-1})$ is an $(r-1)$-chain.  Thus $xw_1\in B_r$, so $\ol{xv_1} = \ol{sv_1'}$ and the map $\widetilde{d}_r$ is well-defined.  
		\end{proof}
		
		\begin{lemma}
			We have $B_r \subset Z_{r+1}= Ker(\widetilde{d}_r)$.
		\end{lemma}
		\begin{proof}
			We only need to show that $Z_{r+1}= Ker(\widetilde{d}_r)$ (the inclusion $B_r\subset Z_{r+1}$ was shown in Lemma \ref{lem:prop_Z_r}). 
			
			For $r=0$ the statement is obvious, so we may assume that $r>0$.
			
			Given $v_{r+1} \in Z_{r+1}$ with a corresponding $(r+1)$-chain $(v_1, \ldots, v_{r+1})$, the truncated sequence $(v_2, \ldots, v_{r+1})$ is an $r$-chain and so $\widetilde{d}_r(v_{r+1})=\ol{xv_2}$. But $xv_2=yv_3$ so $\ol{xv_2}=0$, and so $v_{r+1} \in Ker(\widetilde{d}_r)$. Thus $ Z_{r+1}\subset Ker(\widetilde{d}_r)$.
			
			On the other hand, given $v_{r} \in Ker(\widetilde{d}_r)$ with a corresponding $r$-chain $(v_1, \ldots, v_{r-1}, v_r)$, we may write $xv_1 = yv + xw_1$ where $v\in M$ and $(w_1, \ldots, w_{r-1})$ is an $(r-1)$-chain. Then set $v_i' := v_i - w_i$. We have: $xv_1'=yv$, $yv_i' = yv_i-yw_i = xv_{i+1}-xw_{i+1} = xv_{i+1}'$ for $1\leq i\leq r-2$ and $yv_{r-1}' = yv_{r-1}=xv_r$. Thus $(v, v_1', \ldots, v_{r-1}', v_r)$ is an $(r+1)$-chain and so $v_r \in Z_{r+1}$.
			
			Hence  $ Ker(\widetilde{d}_r) \subset Z_{r+1}$ and the statement is proved.
	\end{proof}}
	
	Thus the map $\widetilde{d}_r$ factors through $E_r = Z_r/B_r$ and we obtain a map 
	\[
	d_r: E_r \longrightarrow E_r.
	\]
	It is not hard to see that $d_r^2=0$ (since $x(xv_1)=0$ already), and so we have obtained our differential.
	
	\InnaA{
		\begin{lemma}
			The cohomology of $d_r$ on $E_{r}$ is isomorphic to $E_{r+1}$.
		\end{lemma}
		\begin{proof}
			We have already seen that $Ker(\widetilde{d}_r) = Z_{r+1}$, so $Ker(d_r) = Z_{r+1}/B_r$.
			
			We now show that $\im(\widetilde{d}_r) = (B_{r+1} + B_r)/B_r$ (the image of $B_{r+1}$ under the quotient map $Z_r \to E_r$). This would prove the required statement, since 
			
			\[\quotient{Z_{r+1}/B_r}{(B_{r+1} + B_r)/B_r} \cong \quotient{Z_{r+1}}{B_{r+1}} = E_{r+1}.\] 
			
			Indeed, let $v_r \in Z_r$ with a corresponding $r$-chain $(v_1, \ldots, v_{r-1}, v_r)$. Then $xv_1 \in B_{r+1}$ and thus $\widetilde{d}_r(v_r) \in (B_{r+1} + B_r)/B_r$. Hence $\im(\widetilde{d}_r) \subset (B_{r+1} + B_r)/B_r$.
			
			Vice versa, any equivalence class in the quotient $(B_{r+1} + B_r)/B_r$ is of the form $\ol{xw_1}$ for some $xw_1 \in B_{r+1}$ with a corresponding $r$-chain $(w_1,\dots,w_r)$. Then $w_r \in Z_r$ and $\widetilde{d}(w_r) = \ol{xw_1}$. Hence $ (B_{r+1} + B_r)/B_r \subset \im(\widetilde{d}_r)$.
		\end{proof}
	}
	\subsection{Leibniz property of \texorpdfstring{$d_r$}{dr}}\label{ssec:spectral_seq_Leibnitz}
	
	We seek to show that the functor of taking the $r$-th page of the above spectral sequence defines a symmetric monoidal functor, for any $r$.  Let $M,N$ be $PGL(1|1)$-modules.  
	
	\begin{lemma}\label{lem:tens_aux}
		We have natural inclusions
		\[
		Z_r(M)\otimes Z_r(N)\sub Z_r(M\otimes N), \ \ \ \ \ B_r(M)\otimes Z_r(N)+Z_r(M)\otimes B_r(N)\sub B_r(M\otimes N).
		\]
	\end{lemma} 
	
	\begin{proof}
		For $r=0,1$ this is clear, so we assume $r>1$. Suppose $v_r\in Z_r(M)$, $w_r\in Z_r(N)$ and \InnaA{let $(v_1, \ldots v_r)$, $(w_1, \ldots, w_r)$ be the corresponding $r$-chains in $M$, $N$ respectively. Then $\left(	\sum_{i+j=k} v_{r-i}\otimes w_{r-j}\right)_{k=r-1, \ldots, 0}$ is an $r$-chain in $M\otimes N$. Indeed,} we have
		\begin{eqnarray*}
			y(v_r\otimes w_r)& = &0 \\ 
			x(v_r\otimes w_r)& = &y(v_{r-1}\otimes w_r+v_r\otimes w_{r-1})\\
			x(v_{r-1}\otimes w_r+v_r\otimes w_{r-1})& = &y(v_{r-2}\otimes w_r+v_{r-1}\otimes w_{r-1}+v_r\otimes w_{r-2})\\
			& \vdots &\\
			\sum_{i+j=k} x(v_{r-i}\otimes w_{r-j}) & = &\sum_{i+j=k+1} y(v_{r-i}\otimes w_{r-j}).
		\end{eqnarray*}
		Thus $Z_r(M)\otimes Z_r(N)\sub Z_r(M\otimes N)$.  Next, we check that 
		\[
		B_r(M)\otimes Z_r(N)+Z_r(M)\otimes B_r(N)\sub B_r(M\otimes N).
		\]
		First we clearly have $\im y\otimes Z_r(N)+Z_r(M)\otimes \im y\sub\im y\rvert_{M\otimes N}$ for any $r\geq 0$.  Next, let \InnaA{$(v_1, \ldots v_{r-1})$ be an $(r-1)$-chain in $M$, and $(w_1, \ldots, w_r)$ be an $r$-chain in $N$}, so that $xv_1 \in B_r(M), w_r \in Z_r(N)$.  We want to show that $xv_1\otimes w_r\in B_r(M\otimes N)$. In the following we work modulo $\im y\rvert_{M\otimes N}$, since it lies in $B_r(M\otimes N)$.
		\begin{eqnarray*}
			xv_1\otimes w_r& = &x(v_1\otimes w_r)-(-1)^{\ol{v_1}}v_1\otimes xw_r\\
			& = &x(v_1\otimes w_r)-(-1)^{\ol{v_1}}v_1\otimes yw_{r-1}\\
			& = &x(v_1\otimes w_r)+yv_1\otimes w_{r-1}\\
			& = &x(v_1\otimes w_r)+xv_2\otimes w_{r-1}\\
			& = &x(v_1\otimes w_r+v_2\otimes w_{r-1})+(-1)^{\ol{v_2}}v_2\otimes xw_{r-1}\\
			& \vdots &\\
			& = &x\left(\sum\limits_{i+j=r-2} v_{r-1-i}\otimes w_{r-j}\right)\pm v_{r-1}\otimes xw_2\\
		\end{eqnarray*}
		
		Now we use that $xw_2=yw_1$ and that $yv_{r-1}=0$ so that the second term $v_{r-1}\otimes xw_2$ lies in $\im y$. Thus we only need to show the first term lies in $B_r(M\otimes N)$. \InnaA{Consider the sequence $\left(\sum\limits_{i+j=k} v_{r-1-i}\otimes w_{r-j} \right)_{k=r-2, r-3, \ldots, 1}$. 
			
			Let us show that this is an $(r-1)$-chain:}
		\begin{eqnarray*}
			y(v_{r-1}\otimes w_{r})& = &0\\
			y(v_{r-2}\otimes w_{r}+v_{r-1}\otimes w_{r-1}) & = &x(v_{r-1}\otimes w_r)\\
			& \vdots & \\
			y(\sum\limits_{i+j=k}v_{r-1-i}\otimes w_{r-j}) & = &x(\sum\limits_{i+j=k-1}v_{r-1-i}\otimes w_{r-j})
		\end{eqnarray*}
		From the above equations we learn that
		\[
		x\left(\sum\limits_{i+j=r-1}v_{r-1-i}\otimes w_{r-j}\right)\in B_r(M\otimes N).
		\]
		Therefore
		\[
		B_r(M)\otimes Z_r(N)\sub B_r(M\otimes N).
		\]
		A similar argument shows that
		\[
		Z_r(M)\otimes B_r(N)\sub B_r(M\otimes N).
		\]
	\end{proof}
	It follows that we have a natural map
	\[
	\Phi_r:E_r(M)\otimes E_r(N)\to E_r(M\otimes N).
	\]
	\begin{prop}
		For each $r\geq 1 $, $\Phi_r$ is an isomorphism.
		
		The action of $d_r$ on $E_r(M\otimes N)$ is by the Leibniz rule:
		\[
		d_r(v_r\otimes w_r)=d_r(v_r)\otimes w_r+(-1)^{\ol{v_r}}v_r\otimes d_r(w_r).
		\]
	\end{prop}
	\begin{proof}

		First of all, recall that given odd operators $d_M: M \to M, d_N: N\to N$ on two supervector spaces $M$ and $N$, we may consider an odd operator $d: M\otimes N \to M\otimes N$ given by the Leibnitz rule $$d(m\otimes n) = d_M(m)\otimes n + (-1)^{\bar{m}} m \otimes d_N(n).$$
		
		Then clearly $Ker(d_M) \otimes Ker(d_N) \subset Ker(d)$ and it is a well-known fact that this induces an isomorphism \begin{equation}\label{eq:ds_tensor}
			\quotient{Ker(d_M)}{Im(d_M)} \otimes \quotient{Ker(d_N)}{Im(d_N)} \cong \quotient{Ker(d)}{Im(d)}
		\end{equation}
		(this is just the statement that $DS$ is a monoidal functor, and follows directly from the representation theory of a $(0|1)$-dimensional Lie superalgebra).
		
		We now prove the statement inductively. 
		
		We know that $d_0=y$ \InnaA{so that for $r=1$, the first statement holds by \eqref{eq:ds_tensor}}. 
		
		Assume that for some $r\geq 1$, $\Phi_r$ an isomorphism. If we can show that $d_r$ satisfies the Leibniz rule, then \InnaA{by \eqref{eq:ds_tensor}, we would conclude that $\Phi_r$ is again an isomorphism}.
		
		\InnaA{Let $(v_1, \ldots v_r)$ be an $r$-chain in $M$, and $(w_1, \ldots, w_r)$ be an $r$-chain in $N$.} We have
		\begin{eqnarray*}
			d_r(v_r\otimes w_r)& = &x(v_1\otimes w_r+v_2\otimes w_{r-1}+\dots+v_r\otimes w_1)\\
			& = &(xv_1\otimes w_r+(-1)^{\ol{v_r}}v_r\otimes xw_1)\\
			& + &[(-1)^{\ol{v_1}}v_1\otimes xw_r+xv_r\otimes w_1+x(v_2\otimes w_{r-1}+\dots+v_{r-1}\otimes w_2)].
		\end{eqnarray*}
		The first term in the above sum is $d_r(v_r)\otimes w_r+(-1)^{\ol{v_r}}v_r\otimes d_r(w_r)$, so it remains to show the second term lies in $B_{r}(M\otimes N)$.  We in fact show that it lies in $\im y$: it is equal to
		\[
		(-1)^{\ol{v_1}}v_1\otimes xw_r+xv_r\otimes w_1+\sum\limits_{i=2}^{r-1}(xv_i\otimes w_{r-i+1}+(-1)^{\ol{v_i}}v_i\otimes xw_{r-i+1})
		\]
		which is equal to:
		\[
		(-1)^{\ol{v_1}}v_1\otimes yw_{r-1}+yv_{r-1}\otimes w_1+\sum\limits_{i=2}^{r-1}(yv_{i-1}\otimes w_{r-i+1}+(-1)^{\ol{v_i}}v_i\otimes yw_{r-i})
		\]
		And this is simply:
		\[
		y\left(v_1\otimes w_{r-1}+v_2\otimes w_{r-2}+\dots+v_{r-1}\otimes w_1\right).
		\]
	\end{proof}
	
	\subsection{Duality}\label{ssec:spectral_seq_duality}  For an odd differential $d$ on a super vector space $M$, let $d^*$ denote the differential on $M^*$ given by
	\[
	d^*(\varphi)(v)=-(-1)^{\ol{\varphi}}\varphi(d(v)).
	\]
	\begin{prop}\label{prop DS dual}
		For each $r\geq0$ we have a natural isomorphism of complexes
		\[
		\Phi_r:(E_r(M^*),d_r(M^*))\to (E_r(M)^*,d_r(M)^*).
		\]
	\end{prop}
	\begin{proof}
		\InnaA{Let $\mathbf{D}\SVec$ be the category of pairs $(V, d)$ where $V\in \SVec$ and $d: V\to \Pi V$ is an odd differential satisfying $d^2= 0$. The maps in this category are maps between vector superspaces (parity preserving) which commute with the differentials.
			
			This category $\mathbf{D}\SVec$ is clearly a rigid (symmetric) monoidal category: these structures are induced by the corresponding structures in $\SVec$.
			
			We now have a functor $E_r: Rep(\mathfrak{pgl}(1|1)) \to \mathbf{D}\SVec$ and we just proved that it is symmetric monoidal.

			By a general argument for monoidal functors (see e.g. \cite{EGNO}, 2.10.6), we obtain a natural isomorphism $E_r(M)^* \cong E_r(M^*)$ for any $M \in Rep(\mathfrak{pgl}(1|1))$.}
	\end{proof}
	
	\subsection{Contragredient duality}\label{sec_appendix_contra_duality} Recall $\sigma_{\g\l(1|1)}$ and $\sigma_n$ from Section \ref{ssec:contra duality}.  We now prove Lemma \ref{lemma twisting DS}:
	
	\begin{lemma}
		We have a natural isomorphism of functors
		\[
		DS^i_{x,y}\circ(-)^{\sigma_{\g\l(1|1)}}\cong (-)^{\sigma_n}\circ DS^i_{y,x}.
		\]
	\end{lemma} 
	\begin{proof}
		For a module $M$, observe that there is a canonical map
		\[
		(DS^i_{y,x}M)^{\sigma_h}\to DS_{x,y}^i(M^\sigma_{\g\l(1|1)})
		\]
		given by $v_r\mapsto v_r$, and it is obviously an isomorphism. Checking that the action of $h$ agrees with this map is straightforward, so we get the result.		
	\end{proof}
	
	\subsection{Action of spectral sequence on finite-dimensional modules}\label{sec_appendix_proof_fds}
	
	We give here the proof of Lemma \ref{action_on_fds}.
	\begin{proof}
		The spectral sequence converges after the first page for all cases but  $DS_{y,x}^r$ acting on $X(n)$ and $DS_{x,y}^r$ acting on $Y(n)$, so these are cases we need to study.  By contragredient duality it's enough to consider $DS_{y,x}^r$ acting on $X(n)$.
		
		For $X(n)$ we have $\im y\sub Z_r(X(n))\sub \ker y$ for all $r>0$, and $\im y$ has codimension 2 in $\ker y$, with complement spanned by $u_{-n+1/2},u_{n-1/2}$, where $u_0$ is even of weight 0 and $u_{n-1/2}$ is odd of weight $n-1/2$.  It is clear that $u_{n-1/2}\in Z_r$ for all $r$.  On the other hand, it is not difficult to check that $u_{-n+1/2}\in Z_r$ only for $r\leq n$.  Thus we have 
		\[
		Z_1=Z_2=\dots=Z_n=\im y+\langle u_{-n+1/2},u_{n-1/2}\rangle, \ \ \ \ Z_r=\im y+\langle u_{n-1/2}\rangle \ \text{ for }r>n.
		\]  
		On the other hand $\im y\sub B_r$ for $r>0$, and $u_{-n+1/2}\notin B_r$ for $r>0$ since $u_{-n+1/2}$ is not in the image of $x$.  However $u_{n-1/2}\in B_{r}$ exactly when $r>n$.  Thus we have
		\[
		E_1=\dots=E_n=\langle u_{-n+1/2},u_{n-1/2}\rangle, \ \ \ \ E_{r}=0 \ \text{ for }r>n.
		\]
		From our setup one can now compute the maps $d_r$, and we find that on $X(n)$ we have: $d_0=y$, $d_{n}:E_n\to E_n$ is given by $d_n(u_{-n+1/2})=u_{n-1/2}$, and $d_r=0$ otherwise.
	\end{proof} 
	
	\bibliographystyle{amsalpha}

	\textsc{\footnotesize Inna Entova-Aizenbud, Department of Mathematics, Ben Gurion University of the Negev, Beer-Sheva,	Israel} 
	
	\textit{\footnotesize Email address:} \texttt{\footnotesize entova@bgu.ac.il}
	
	\textsc{\footnotesize Vera Serganova, Dept. of Mathematics, University of California at Berkeley, Berkeley, CA 94720} 
	
	\textit{\footnotesize Email address:} \texttt{\footnotesize serganov@math.berkeley.edu}
	
	\textsc{\footnotesize Alexander Sherman, School of Mathematics and Statistics, University of Sydney, Camperdown NSW 2006} 
	
	\textit{\footnotesize Email address:} \texttt{\footnotesize xandersherm@gmail.com}

\end{document}